\documentclass[12pt]{amsart}

\usepackage{amssymb}
\usepackage{latexsym}
\usepackage[enableskew]{youngtab}
\usepackage{verbatim}
\usepackage{fullpage}
\usepackage{pslatex}

\input {cyracc.def}
\font\ququ=cmr10 scaled \magstep1
\font\tencyr=wncyr10 
\font\tencyi=wncyi10 
\font\tencysc=wncysc10 
\def\rus{\tencyr\cyracc}
\def\rusi{\tencyi\cyracc}
\def\rusc{\tencysc\cyracc}

\newcommand{\re}[1]{\textrm  (\ref{#1})}

\renewenvironment{proof}
{\noindent {\sl Proof.}\quad }{\hfill
$\square$ \vskip1.1ex\noindent }

\newenvironment{proof*}
{\noindent {\sl Proof.}\quad }{\hfill
$\square$}

\renewcommand{\theequation}{\thesection .\arabic{equation}}
\renewcommand{\thesubsubsection}{\theequation .\arabic{subsubsection}}

\catcode`\@=\active
\catcode`\@=11
\def\@eqnnum{\hbox to
.01pt{}\rlap{\hskip-\displaywidth(\mathbf{\theequation})}}
\catcode`\@=12

\newenvironment{s}[1]
{ \vskip1.2ex \refstepcounter{equation}
\noindent {\bf \theequation\enspace #1.} \begin{sl}}{\end{sl}
\vskip1.1ex\noindent }

\newenvironment{rem}[1]
{ \vskip1.2ex \refstepcounter{equation}
\noindent {\bf \theequation\enspace {#1}.} }{ \vskip1.1ex\noindent }

\newcommand {\ah}{{\frak a}}
\newcommand {\be}{{\frak b}}
\newcommand {\ce}{{\frak c}}

\newcommand {\ee}{{\frak e}}

\newcommand {\g}{{\frak g}}

\newcommand {\fH}{{\goth H}}

\newcommand {\te}{{\frak t}}
\newcommand {\ut}{{\frak u}}

\newcommand {\slno}{{\frak sl}_{n+1}}

\newcommand {\spn}{{\frak sp}_{2n}}
\newcommand {\sono}{{\frak so}_{2n+1}}
\newcommand {\sone}{{\frak so}_{2n}}

\newcommand {\gA}{{\goth A}}
\newcommand {\gC}{{\goth C}}
\newcommand {\gH}{{\goth H}}

\newcommand {\esi}{\varepsilon}
\newcommand {\ap}{\alpha}

\newcommand {\lb}{\lambda}
\newcommand {\vp}{\varphi}

\newcommand {\HW}{\widehat W}
\newcommand {\HV}{\widehat V}
\newcommand {\HP}{\widehat\Pi}
\newcommand {\HD}{\widehat\Delta}

\newcommand {\EE}{{\mathcal E}}

\newcommand {\ck}{{\mathcal K}}
\newcommand {\MM}{{\mathcal M}}
\newcommand {\PP}{{\mathcal P}}
\newcommand {\N}{{\mathcal N}}

\newcommand {\VV}{{\Bbb V}}

\newcommand {\ad}{{\mathrm{ad\,}}}
\newcommand {\Ad}{{\mathrm{Ad\,}}}
\newcommand {\Aut}{{\mathrm{Aut\,}}}

\newcommand {\rk}{{\mathrm{rk\,}}}

\newcommand {\tri}{{\frak sl}_2}
\newcommand {\GR}[2]{{\textrm{{\bf #1}}}_{#2}}

\newcommand {\Ab}{{\frak Ab}}
\newcommand {\AD}{{\frak Ad}}

\newcommand {\EAD}{{\mathcal E}(\AD)}
\newcommand {\EAb}{{\mathcal E}(\Ab)}
\newcommand {\EADg}{{\mathcal E}(\AD(\g))}
\newcommand {\EAbg}{{\mathcal E}(\Ab(\g))}

\newcommand {\beq}{\begin{equation}}
\newcommand {\eeq}{\end{equation}}

\newcommand{\curge}{\succcurlyeq}
\newcommand{\curle}{\preccurlyeq}
\renewcommand{\le}{\leqslant}
\renewcommand{\ge}{\geqslant}

\newcommand{\vts}{\VV_{\theta_s}}

\newcommand{\adn}{{\sf ad}-nilpotent }
\newcommand{\apo}{\mbox{$\ap_0$}}
\newcommand{\apx}{\mbox{$\ap_1$}}
\newcommand{\apxx}{\mbox{$\ap_2$}}
\newcommand{\apxxx}{\mbox{$\ap_3$}}
\newcommand{\apxxxx}{\mbox{$\ap_4$}}
\newcommand{\apxxxxx}{\mbox{$\ap_5$}}
\newcommand{\apxxxxxx}{\mbox{$\ap_6$}}
\newcommand{\zam}{\mbox{$\gamma_{11}$}}
\newcommand{\zan}{\mbox{$\gamma_{12}$}}
\newcommand{\zap}{\mbox{$\gamma_{13}$}}
\newcommand{\zaq}{\mbox{$\gamma_{14}$}}
\newcommand{\zar}{\mbox{$\gamma_{15}$}}
\newcommand{\zas}{\mbox{$\gamma_{16}$}}
\newcommand{\zbn}{\mbox{$\gamma_{22}$}}
\newcommand{\zbp}{\mbox{$\gamma_{23}$}}
\newcommand{\zbq}{\mbox{$\gamma_{24}$}}
\newcommand{\zbr}{\mbox{$\gamma_{25}$}}
\newcommand{\zbs}{\mbox{$\gamma_{26}$}}
\newcommand{\zcp}{\mbox{$\gamma_{33}$}}
\newcommand{\zcq}{\mbox{$\gamma_{34}$}}
\newcommand{\zcr}{\mbox{$\gamma_{35}$}}
\newcommand{\zcs}{\mbox{$\gamma_{36}$}}
\newcommand{\zdq}{\mbox{$\gamma_{44}$}}
\newcommand{\zdr}{\mbox{$\gamma_{45}$}}
\newcommand{\zds}{\mbox{$\gamma_{46}$}}
\newcommand{\zer}{\mbox{$\gamma_{55}$}}
\newcommand{\zes}{\mbox{$\gamma_{56}$}}
\newcommand{\zfs}{\mbox{$\gamma_{66}$}}

\newfam\Bbbfam\newfam\eufam\newfam\eusfam%
\font\Bbbfont=msbm10 scaled 1200%
\font\olala=msam10 scaled 1200%
\font\Bbbsmallfont=msbm8%
\textfont\Bbbfam=\Bbbfont\scriptfont\Bbbfam=\Bbbsmallfont
\font\euzw=eufm10 scaled 1200%
\font\euac=eufm7 scaled 1200%
\font\euacc=eufm7 scaled 1000%
\textfont\eufam=\euzw\scriptfont\eufam=\euac%
\scriptscriptfont\eufam=\euacc%
\font\euszw=eusm10 scaled 1200%
\font\eusac=eusm7 scaled 1200%
\font\eusacc=eusm7 scaled 1000%
\textfont\eusfam=\euszw\scriptfont\eusfam=\eusac%
\scriptscriptfont\eusfam=\eusacc%
\def\frak{\fam\eufam}%
\def\goth{\fam\eusfam}%
\def\Bbb{\fam\Bbbfam}%

\def\varnothing{\hbox {\Bbbfont\char'077}}
\def\square{\hbox {\olala\char"03}}

\begin{document}
\setlength{\parskip}{2pt plus 4pt minus 0pt}
\hfill {\scriptsize November 24, 2004} 
\vskip1ex
\vskip1ex

\title[The poset of positive roots]{The poset of positive roots and its relatives}
\author{\sc Dmitri I. Panyushev}
\thanks{This research was supported in part by 
CRDF Grant no. RM1-2543-MO-03}
\keywords{Simple Lie algebra, root system, Hasse diagram, ad-nilpotent ideal}
\subjclass[2000]{17B20, 20F55}
\maketitle
\begin{center}
{\footnotesize
{\it Independent University of Moscow,
Bol'shoi Vlasevskii per. 11 \\
119002 Moscow, \quad Russia \\ e-mail}: {\tt panyush@mccme.ru }\\
}

\vskip2ex
{\bf Abstract}
 
{\small
\parbox{410pt}{
Let $\Delta$ be a root system with a subset of positive roots, 
$\Delta^+$. We consider edges of the Hasse diagrams of some
posets associated with $\Delta^+$. For each edge one naturally 
defines its {\it type\/}, and we study the partition of the
set of edges into types. For $\Delta^+$, the type is a simple root,
and for the posets of {\sf ad}-nilpotent and Abelian ideals the type
is an affine simple roots. We give several descriptions of the
set of edges of given type and uniform expressions for the number of edges. 
By a result of Peterson, the number of Abelian ideals is $2^n$, where $n$ is the
rank of $\Delta$. We prove that the number of edges of the corresponding
Hasse diagram is $(n+1)2^{n-2}$.
For $\Delta^+$ and the Abelian ideals, we compute the number of edges of each type
and prove that the number of edges of type $\ap$ depends only on
the length of the root $\ap$.
}
}
\end{center}

\vskip2ex
\noindent
Let $\Delta$ be a reduced irreducible root system in an $n$-dimensional real
vector space $V$. Choose a subsystem of positive roots $\Delta^+$ with the corresponding
subset of simple roots $\Pi=\{\ap_1,\ldots,\ap_n\}$.
Let $G$ be the corresponding
simply-connected simple algebraic group with Lie algebra $\g$.
Fix a triangular decomposition $\g=\ut^+ \oplus\te\oplus\ut^-$, where $\te$
is a Cartan subalgebra and $\ut^+=\underset{\gamma\in\Delta^+}{\bigoplus}\g_\gamma$.
Then $\be=\te\oplus\ut^+$ is the fixed Borel subalgebra.
The root order in $V$ is given  by letting $x\curle y$ if
$y-x$ is a non-negative integral combination of positive roots. 
In particular, we always regard
$\Delta^+$ as  poset under `$\curle$'.
Then the highest root, $\theta$, is the unique maximal
element of $\Delta^+$, whereas the simple roots are precisely the minimal elements.

In this article, we consider combinatorial properties of several posets associated
with $\Delta^+$. We will especially be interested in the edges of the
Hasse diagrams of these posets. 
Given a poset $\PP$, we write $\gH(\PP)$
for the Hasse diagram of $\PP$ and $\mathcal E(\PP)$ for the set of
edges of $\gH(\PP)$. 
The common property of all posets to be considered below
is that one can naturally define the {\it type\/} of any edge of $\gH(\PP)$.
If $T$ is the parameter set of possible types, then
$\mathcal E(\PP)=\underset{t\in T}{\sqcup} \mathcal E(\PP)_t$, and our aim is to 
describe this partition.

In the simplest case, the poset in question is $\Delta^+$ itself.
If $\nu$ covers $\mu$ in $\Delta^+$, then the type of the edge $(\nu,\mu)$
is the simple root $\nu-\mu$.
Let $h$ (resp. $h^*$) denote the Coxeter (resp. dual Coxeter) number
of $\Delta$.
We show that if $\ap$ is long, then $\#\mathcal E{(\Delta^+)}_\ap=h^*-2$;
if $\ap$ is short, then  $\#\mathcal E(\Delta^+)_\ap$
is equal to either to $h-2$ or $h-3$, depending on $\Delta$.
Anyway, $\#\mathcal E(\Delta^+)_\ap$  depends only on $\|\ap\|$, the length
of $\ap$. Although these results are not difficult, they are
apparently new, see Section~\ref{odin}.
A similar phenomenon occurs in several related cases. First, if $\Delta^+$
is not simply-laced, then
one can consider the subposet of $\Delta^+$ consisting of all short roots.
More generally, given a simple $G$-module $\VV$, one can consider the
poset of weights of $\VV$, see Section~\ref{dva}.

Another interesting poset is that of all \adn ideals in $\be$, 
denoted by $\AD$. There have recently been a lot of activity in
studying \adn ideals \cite{cp1,cp2,cp3,ko1,lp,duality, mima-p,eric,suter}.
Combinatorially, $\AD$ is the set of all {\it upper\/}
(\,=\,dual order) ideals of $\Delta^+$, 
the partial order being given by the usual containment.
An edge of $\gH(\AD)$ is a pair of \adn ideals $(\ce,\ce')$
such that $\ce'\subset\ce$ and $\dim\ce=\dim\ce'+1$. 
We then say that the edge $(\ce,\ce')$ {\it terminates\/} in $\ce$.
It is easily seen that the edges terminating in $\ce$ are in a bijection with
the generators of $\ce$, i.e., the minimal roots in the corresponding 
upper ideal $I_\ce$.
To define the type of edges in this situation, we exploit the following
ingredients:

$\bullet$ \ a bijection
between the \adn ideals and certain elements of the affine Weyl group, $\HW$
({\sc Cellini--Papi}~\cite{cp1}; in the Abelian case this goes back to {\sc Peterson}
and {\sc Kostant}~\cite{ko1}). Given $\ce\in\AD$, let $w_\ce$ denote the corresponding
element of $\HW$. 

$\bullet$ \ the characterisation of the generators of
$\ce$ in terms of $w_\ce$ \cite{duality}.
\\[.6ex]
Now, the type of an edge appears to be an {\it affine\/} simple root.
The set of affine simple roots is
$\HP=\{\ap_0\}\cup\Pi$, where $\ap_0$ is an extra root, see details in
Section~\ref{pos_all}.
We give a closed formula for the number of edges and provide several geometric
descriptions for the set of edges of a given type.
The first of them is based on a relationship between $\AD$ and the integral points 
of the simplex $D_{min}$; second description uses the action of $w_\ce$ on the 
fundamental alcove in $V$, see Theorems~\ref{tip-edges} and \ref{geom}.
The last description is given in terms of reduced decompositions
of $w_\ce$.  Namely,  $\fH(\AD)$ has an edge of type
$\ap_i$ terminating in $\ce$
if and only if $w_\ce$ has a reduced decomposition starting with $s_i$.
Here $s_i\in\HW$ is the reflection corresponding to $\ap_i\in\HP$.

The subposet of $\AD$ consisting of all Abelian ideals is denoted by $\Ab$.
Because $\fH(\Ab)$ is a subdiagram of $\fH(\AD)$, the type of edges
in $\fH(\Ab)$ is well-defined.
By a result of D.\,Peterson, $\#\Ab=2^n$ (see \cite{cp1},\cite{ko1}). 
We first show conceptually that the number of edges of 
$\fH(\Ab)$ equals $(n+1)2^{n-2}$.
Our next result asserts that $\#\mathcal E(\Ab)_\ap$, $\ap\in\HP$,
depends only on $\|\ap\|$.
Furthermore, if $\g\ne\spn$, then the number of edges of each type
is $2^{n-2}$. Unfortunately, the proof consists of explicit verifications 
for all simple Lie algebras. To this end, we choose the following path. 
A root $\gamma\in\Delta^+$ is said to be {\it commutative\/} if the upper ideal
generated by $\gamma$ is Abelian. The commutative roots form an upper
ideal, which is not necessarily Abelian. (We notice that the complement of this
ideal has a unique maximal element, which can explicitly be described; a uniform 
description of the cardinality of the set of commutative roots is also given,
see Theorem~\ref{number:com}.) The utility of commutative roots is revealed via
the following property. 
Suppose $\ah\in\Ab$ and $\gamma$ is a generator of $\ah$. 
By the general rule mentioned above for $\AD$, the pair $(\ah,\gamma)$ determines
an edge terminating in $\ah$. If this edge is of type $\ap$, then
we say that $\gamma$ has {\it class\/} $\ap$ (in the Abelian ideal $\ah$).
The point here is that the class of $\gamma$ does not depend on an Abelian
ideal in which $\gamma$ is a generator. Thus, each commutative root gains a
well-defined class, which is an element of $\HP$.
Furthermore, we provide effective methods
of computing the class of commutative roots. Our strategy for computing
the numbers $\#\mathcal E(\Ab)_\ap$, $\ap\in\HP$,
begins with determining all commutative
roots of class $\ap$, say $\{\gamma_1,\ldots,\gamma_t\}$. Then we compute, 
for each $i$,
the number of Abelian ideals having $\gamma_i$ as a generator.
The sum of all these numbers equals $\#\mathcal E(\Ab)_\ap$.
For the classical Lie algebras the description of classes of the commutative roots
is given using the usual matrix realisations; for the exceptional cases, we list
the classes in the Appendix.

In Section~\ref{fin_poset}, we introduce the covering polynomial of a
finite poset, and compute this polynomial for the posets considered in this article.

{\small
{\bf Acknowledgements.} Part of this work was done while I was visiting
the Max-Planck-Institut f\"ur Mathematik (Bonn). I thank the Institute for
its hospitality and inspiring environment.
}


\section{Edges of the Hasse diagram of $\Delta^+$}
\label{odin}
\setcounter{equation}{0}

\noindent
Let $\fH={\fH}(\Delta^+)$ be the Hasse diagram of $(\Delta^+, \curle)$.
In other words, $\fH$ is a directed graph whose set of vertices
is $\Delta^+$, and the set of edges, ${\mathcal E}(\fH)$, consists of the
pairs of positive roots $(\mu,\nu)$ such that $\mu-\nu\in\Pi$.
We then say that $(\mu,\nu)$ is an {\it edge of type\/} $\mu{-}\nu$.

Let $h=h(\Delta)=h(\g)$ denote the Coxeter number of $\Delta$ (or $\g$). As is well-known,
the number of vertices of $\fH$ equals $nh/2$.

\begin{s}{Theorem}  \label{simply}
Suppose $\Delta$ is simply-laced. Then the number of edges of each type
in ${\fH}(\Delta^+)$ is equal to $h-2$. The total number of edges 
equals $n(h-2)$.
\end{s}\begin{proof}
Clearly, the number of edges of type $\ap$ is 
equal to the number of positive roots $\nu$ such that $\nu+\ap\in\Delta^+$.
Since all roots have the same length, the latter is the same
as the number of positive roots $\nu$ such that $(\nu,\ap)<0$.
\\
For any $\gamma\in\Delta$, we set $\Delta(\gamma)=
\{\mu\in\Delta \mid (\mu,\gamma)\ne 0\}$.
By \cite[chap.\,VI,~\S\,11,~Prop.\,32]{bour}, we have $\#\Delta(\gamma)=4h-6$ for any $\gamma$.
Therefore the number of such positive roots is $2h-3$.
Consider the partition of this set according to the sign of the scalar product
with $\gamma$:
\[
   \Delta(\gamma)^+=\Delta(\gamma)^+_{>0}\sqcup \Delta(\gamma)^+_{<0} \ .
\]
Now, let $\gamma=\ap$ be a simple root. Then 
$\ap\in \Delta(\ap)^+_{>0}$, and it is easily seen that the reflection
$s_\ap$ yields a bijection:
\[
   s_\ap:  \Delta(\ap)^+_{>0}\setminus\{\ap\} \to \Delta(\ap)^+_{<0} \ .
\]
Hence the $\#\Delta(\ap)^+_{<0}=h-2$. But this is exactly the number of edges of
type $\ap$.
\end{proof}%
If $\Delta$ is not simply-laced, then one has to distinguish long and short roots.
Let $\Pi_l$ (resp. $\Pi_s$) stand for the set of long (resp. short) simple
roots.
Let $h^*=h^*(\Delta)$ be the dual Coxeter number of $\Delta$.
By definition, $h^*=(\rho,\theta^\vee)+1$. If $\Delta$ is simply-laced, then
$h^*=h$. The values of $h$ and $h^*$ in the non-simply-laced cases are given in 
the following table
\vskip.7ex

\begin{center}\begin{tabular}{c|cccc|}
      & $\GR{B}{n}$ & $\GR{C}{n}$ & $\GR{F}{4}$ & $\GR{G}{2}$ \\ \hline
$h$   &   $2n$      &    $2n$     &    12       &     6       \\
$h^*$ &   $2n{-}1$  &    $n{+}1$  &     9       &     4       
\end{tabular}\end{center} 

\begin{s}{Theorem}  \label{nons}
Suppose $\Delta$ is not simply-laced.
\begin{itemize}
\item[\sf (i)] \ 
For $\ap\in\Pi_l$, the number of edges of type $\ap$ in ${\fH}(\Delta^+)$
equals\/ $h^*{-}2$.
\item[\sf (ii)] \ 
For $\ap\in\Pi_s$, the number of edges of type $\ap$ in 
${\fH}(\Delta^+)$ equals  \\
$\left\{ \begin{array}{cl} h-2, & \text{if } \ \Delta\in \{\GR{B}{n},\GR{C}{n},
\GR{F}{4}\}, \\  h-3, &\text{if } \ \Delta=\GR{G}{2} \quad .
\end{array}\right.$ 
\end{itemize}
\end{s}\begin{proof}
(i) If $\ap\in\Pi_l$, then the number of edges of type $\ap$ is again equal
to the number of $\nu\in\Delta^+$ such that $(\nu,\ap)<0$.
Therefore we can argue as in Theorem~\ref{simply}. The only difference is that 
now we refer to \cite{rudi} for the equality $\#\Delta(\ap)=4h^*{-}6$.

(ii) For $\ap\in\Pi_s$, it is no longer true that the number of edges of type $\ap$ 
equals the number of $\nu\in\Delta^+$ such that $(\nu,\ap)<0$. For, it may happen that
$(\nu,\ap)\ge 0$, but $\nu+\ap$ is a root. 
So that we give a case-by-case argument. Our notation and numbering of simple roots
follows \cite[Tables]{VO}.
\\[.5ex]
{\bf --} \ In case of $\GR{B}{n}$, we have $\ap_n=\esi_n$ is the only short simple root.
The set of roots $\mu$ such that $\mu+\ap_n$ is a root consists of
$\esi_i-\esi_n \ (i<n)$, $\esi_1,\ldots, \esi_{n-1}$. Thus, there are $2n{-}2$ 
such roots.
\\[.5ex]
{\bf --} \ In case of $\GR{C}{n}$, we have $\ap_i=\esi_i-\esi_{i+1}\in\Pi_s$, $(i<n)$.
The set of roots $\mu$ such that $\mu+\ap_i$ is a root consists of
$\esi_k-\esi_i \ (k<i)$; $\esi_{i+1}-\esi_l \ (i{+}1<l)$;
$\esi_{i+1}+\esi_{m}$ for all $m$. Thus, there are $2n{-}2$ 
possibilities for each $i$.

The two exceptional cases are left to the reader. (The Hasse diagram for
the $\GR{F}{4}$-case is depicted below.)
\end{proof}%
{\bf Remarks.} 1. In the following section we prove {\it a priori\/}
(in a slightly more general context)
that the number of edges depends only on the length of a simple root, see
Theorem~\ref{VV-graph}.

2. It is not hard to find an isomorphism between the posets $\Delta^+(\GR{B}{n})$ and
$\Delta^+(\GR{C}{n})$. Their Hasse diagrams are therefore isomorphic, too, 
and hence have the same number of edges. However, this isomorphism does not respect 
the length of roots and type of edges.

3. The second claim of Theorem~\ref{nons} for $\Delta\in \{
\GR{B}{n}, \GR{C}{n}, \GR{F}{4}\}$ admits a partial explanation.
Suppose $\ap\in\Pi_s$ and $e_\ap\in \g_\ap$ is a nonzero root vector.
Then $(\ad e_\ap)^3=0$, and an easy calculation shows that the number of edges of 
type $\ap$ is equal to $(\dim G{\cdot}e_\ap-2)/2$.
For these three cases, there is an interesting phenomenon observed by R.~Brylinski and B.~Kostant \cite{BK}.
The orbit $G{\cdot}e_\ap\subset\g$ is {\it shared} in the following sense.
Let $\tilde\g$ be the simple Lie algebra obtained from (the Dynkin diagram of) $\g$
by unfolding, i.e., $\GR{B}{n}\mapsto \GR{D}{n+1}$, $\GR{C}{n} \mapsto 
\GR{A}{2n-1}$, $\GR{F}{4}\mapsto \GR{E}{6}$.
Then the minimal nilpotent
$\tilde G$-orbit in $\tilde\g$ is a finite covering of $G{\cdot}e_\ap$.
Furthermore, we have $h(\g)=h(\tilde \g)$. Finally, since $\tilde\g$ is simply-laced,
the dimension of the 
minimal nilpotent orbit in $\tilde\g$ equals $2h(\tilde\g)-2$.

\begin{figure}[htb]
\setlength{\unitlength}{0.025in}
\begin{center}
\begin{picture}(100,165)(1,30)

\put(50,190){\circle*{2}}
\put(65,175){\circle*{2}}
\put(70,160){\circle*{2}}
\put(65,145){\circle{2}}
\put(50,130){\circle{2}}
\put(60,130){\circle*{2}}
\put(45,115){\circle{2}}     
\put(65,115){\circle*{2}}
\put(30,100){\circle*{2}}     
\put(50,100){\circle{2}}    
\put(80,100){\circle*{2}}
\put(35,85){\circle*{2}}     
\put(45,85){\circle{2}}    
\put(65,85){\circle{2}}
\put(30,70){\circle{2}}     
\put(50,70){\circle*{2}}    
\put(60,70){\circle{2}}
\put(25,55){\circle*{2}}     
\put(45,55){\circle{2}}    
\put(65,55){\circle{2}}
\put(30,40){\circle*{2}}     
\put(40,40){\circle*{2}}    
\put(50,40){\circle{2}}    
\put(60,40){\circle{2}}
                                      \put(60.7,131.4){\line(1,3){4}}
\put(51,131){\line(1,1){12.7}}        \put(45.7,116.4){\line(1,3){4}} 
\put(46,116){\line(1,1){12.7}}        \put(65.7,146.4){\line(1,3){4}}
\put(31,101){\line(1,1){12.7}}        \put(45.7,86.4){\line(1,3){4}}
\put(51,101){\line(1,1){12.7}}        \put(30.7,71.4){\line(1,3){4}} 
\put(36,86){\line(1,1){12.7}}         \put(45.7,56.4){\line(1,3){4}}   
\put(66,86){\line(1,1){12.7}}         \put(25.7,56.4){\line(1,3){4}}    
\put(51,71){\line(1,1){12.7}}         \put(60.7,71.4){\line(1,3){4}}    
\put(51,41){\line(1,1){12.7}}         \put(40.7,41.4){\line(1,3){4}}    
\put(46,56){\line(1,1){12.7}}         \put(60.7,41.4){\line(1,3){4}}    
\put(31,71){\line(1,1){12.7}}           

\put(51,189){\line(1,-1){12.8}}
\put(66,114){\line(1,-1){12.8}}          
\put(51,99){\line(1,-1){12.8}}        \put(69.3,161.3){\line(-1,3){4}} 
\put(36,84){\line(1,-1){12.8}}        \put(64.3,116.3){\line(-1,3){4}}
\put(31,69){\line(1,-1){12.8}}        \put(49.3,101.3){\line(-1,3){4}}
\put(26,54){\line(1,-1){12.8}}        \put(34.3,86.3){\line(-1,3){4}} 
\put(46,84){\line(1,-1){12.8}}        \put(64.3,56.3){\line(-1,3){4}}   
                                      \put(29.3,41.3){\line(-1,3){4}}    
                                      \put(49.3,41.3){\line(-1,3){4}}    
\put(27,33){{\small $\ap_4$}}
\put(37,33){{\small $\ap_3$}}
\put(47,33){{\small $\ap_2$}}
\put(57,33){{\small $\ap_1=[1000]$}}  
\put(29,188){{\small $[2432]$}}    
\put(59,182){{\small $\ap_4$}} \put(69,167){{\small $\ap_3$}} 
\put(69,150){{\small $\ap_2$}} \put(51,138){{\small $\ap_1$}} 
\put(14,83){{\small $[0211]$}}    
\put(9,98){{\small $[0221]$}}    
\put(28,128){{\small $[1321]$}}    
\put(23,113){{\small $[1221]$}}    
\put(70,53){{\small $[1100]$}}   
\put(65,68){{\small $[1110]$}}   
\put(70,83){{\small $[1210]$}}   
\put(84,98){{\small $[2210]$}}   
\put(70,113){{\small $[2211]$}}   
\put(64,128){{\small $[2221]$}}  
\put(69,143){{\small $[2321]$}} 
\put(74,158){{\small $[2421]$}} 

\end{picture}
\end{center}
\caption{The Hasse diagram of $\Delta^+(\GR{F}{4})$}   \label{pikcha_F4}
\end{figure}
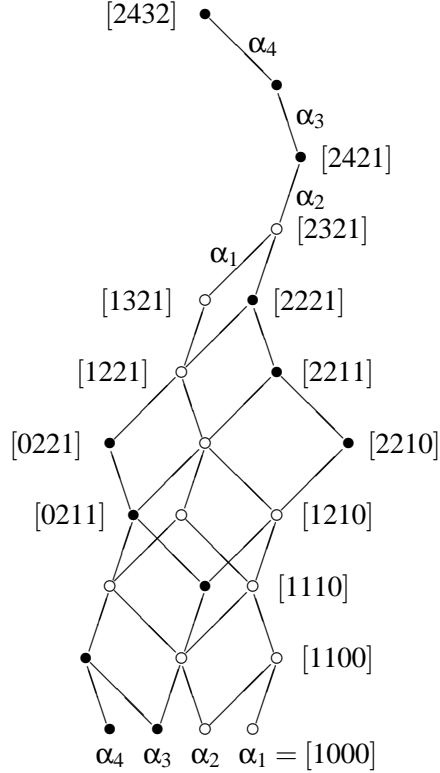

\begin{rem}{Example}
The Hasse diagram $\fH(\Delta^+(\GR{F}{4}))$ has 34 edges,
see Figure~\ref{pikcha_F4}. 
The edges of the same type are drawn to have the same slope. Therefore we have 
indicated only the type of 4 edges in the upper part of the diagram. The black nodes
represent long roots.
For brevity, the root $\sum_{i=1}^4 m_i\ap_i$ is denoted by $[m_1\,m_2\,m_3\,m_4]$.
\end{rem}
\noindent


\section{Some generalisations}
\label{dva}
\setcounter{equation}{0}

\noindent                
Suppose $\Delta^+$ has roots of different length.
We use subscripts
`s` and `l` to mark various objects related to short and long roots, respectively.
For instance, $\Delta_s$ is the set of short roots, $\Delta=\Delta_s\sqcup
\Delta_l$, and $\Pi_s=\Pi\cap \Delta_s$. Set $\rho_s=\frac{1}{2}\vert
\Delta^+_s\vert$ and $\rho_l=\frac{1}{2}\vert
\Delta^+_l\vert$.
Let $\theta\in\Delta^+$ be the highest root and $\theta_s$ the short dominant root.
Recall that $\Delta_l=W{\cdot}\theta$, $\Delta_s=W{\cdot}\theta_s$, and 
${|\!|\theta|\!|^2}/{|\!|\theta_s|\!|^2}=2$ or~3. 
Consider $\Delta^+_s$ as subposet of $\Delta^+$.
This poset has a geometric meaning, see below.
Let me stress that
although the vertices of $\fH(\Delta^+_s)$ represent the short roots only, the edges 
still correspond to the whole of $\Pi$. 
It follows from \cite[ch.VI,\,\S\,1,\,Prop.\,33]{bour} that
$\#(\Delta^+_s)=h{\cdot}\#(\Pi_s)$.
Let us now compute the number of edges of each type.
 
\begin{s}{Theorem}
\begin{itemize}
\item[\sf (i)] \ 
For $\ap\in\Pi_s$, the number of edges of type $\ap$ 
in ${\fH}(\Delta^+_s)$ is\/ $h(\Delta_s){-}2$.
\item[\sf (ii)] \ 
For $\ap\in\Pi_l$, the number of edges of type $\ap$ in ${\fH}(\Delta^+_s)$
is  $h^*(\Delta)- h(\Delta_l)$.
\end{itemize}
\end{s}\begin{proof}
(i) An edge of type $\ap\in\Pi_s$ is a pair $(\gamma,\mu)\in\Delta^+_s\times\Delta^+_s$ such that
$\gamma-\mu=\ap$. Hence the number of edges of type $\ap$ is equal to
$\#\{\mu\in\Delta^+_s\mid (\mu,\ap)<0\}$. Recall that $\Delta_s$ is a root 
system in its own right. It can be reducible, but all irreducible subsystems 
are isomorphic, since $W$ acts transitively on the roots of the same length.
Therefore the Coxeter number $h(\Delta_s)$ is well-defined.
Furthermore, $\Pi_s$ forms a part of a basis for  $\Delta^+_s$.
Hence one may refer to Theorem~\ref{simply}.

(ii) \ If $\ap\in\Pi_l$, then again the number of edges of type $\ap$ 
is equal to 
\begin{multline*}
\#\{\mu\in\Delta^+_s\mid (\ap,\mu)<0 \}=
\#\{\mu\in\Delta^+\mid (\ap,\mu)<0 \}-\#\{\mu\in\Delta^+_l\mid (\ap,\mu)<0\}=
\\
=(h^*(\Delta)-2)-(h(\Delta_l)-2)=h^*(\Delta)-h(\Delta_l) \ .
\end{multline*}
Here we have also used the fact that $\Delta_l$ is a root
system and $h(\Delta_l)$ is well-defined.
\end{proof}%
{\bf Remark.} There is the obvious recipe for obtaining $\fH(\Delta^+_s)$ from
$\fH(\Delta^+)$. One has to erase all vertices (with connecting edges)
corresponding to the long roots.

\vskip1ex
It is easy to realise that one can consider similar problems for arbitrary
representations of simple Lie algebras. 
Let $\VV_\lb$ be the simple finite-dimensional $\g$-module with highest 
weight $\lb$. We write $m_\lb(\mu)$ for the dimension of $\mu$-weight
space of $\VV_\lb$.
The set of weights of $\VV_\lb$ is denoted by ${\mathcal P}(\VV_\lb)$.
Being a subset of $V$, ${\mathcal P}(\VV_\lb)$ can be regarded as poset under
`$\curle$'.
The {\it weight diagram\/} of $\VV_\lb$, 
denoted ${\frak W}(\VV_\lb)$, is the directed multigraph whose set 
of vertices is ${\mathcal P}({\VV_\lb})$ and the set of edges consists
of the pairs $(\nu,\mu)$ such that $\nu-\mu\in\Pi$. The multiplicity
of the edge $(\nu,\mu)$ is defined as \ $\min\{m_\lb(\mu), m_\lb(\nu)\}$.
It is easy to check that ${\frak W}(\VV_\lb)$ is a simple graph
(i.e., all the multiplicities are equal to 1) if and only if
$m_\lb(\mu)=1$ for all nonzero $\mu\in {\mathcal P}({\VV_\lb})$.
If ${\frak W}(\VV_\lb)$ is a simple graph, then it is nothing but the
Hasse diagram of ${\mathcal P}({\VV_\lb})$.
As above, we define the ``type'' of each edge of ${\frak W}(\VV_\lb)$.

\begin{s}{Theorem}   \label{VV-graph}
If $\ap,\beta\in \Pi$ have the same length,
then the number of edges of type $\ap$ and $\beta$ in ${\frak W}(\VV_\lb)$
(counted with multiplicities) is the same.
That is to say, the number of edges of type $\ap\in\Pi$  depends only on
the length of $\ap$. 
\end{s}\begin{proof*}
Let $\tri(\ap)$ be the simple three-dimensional subalgebra of $\g$ 
corresponding to $\ap$. Consider $\VV_\lb$ as $\tri(\ap)$-module.
We write $R(d)$ for the simple $\tri$-module of dimension $d+1$.
It is easy to check that if $\VV_\lb\vert_{\tri(\ap)}\simeq \oplus n_i R(d_i)$,
then the number of edges of type $\ap$ is equal to $\sum_i n_i d_i$.
On the other hand, the subalgebras $\tri(\ap)$ and $\tri(\beta)$ are 
conjugate under $\Aut \g$, so that $\VV_\lb\vert_{\tri(\ap)}$
and $\VV_\lb\vert_{\tri(\mu)}$ are isomorphic as $\tri$-modules.
\end{proof*}%
\begin{rem}{Examples}
1. If $\VV_\lb=\g$ is the adjoint $\g$-module, then 
${\mathcal P}(\g)=\Delta\cup\{0\}$. Here ${\frak W}(\g)$ is a simple graph.
After deleting $\{0\}$ from this weight diagram, we obtain two isomorphic
connected components corresponding to $\Delta^+$ and $\Delta^-$.
This provides another (a priori) proof for the fact that the number of
edges in $\fH(\Delta^+)$ of type $\ap\in\Pi$ depends only on the length
of $\ap$. Clearly, we loose two edges of each type with deleting the
vertex $\{0\}$. Therefore it follows from Theorem~\ref{nons}
that the number of edges in ${\frak W}(\g)$ of type $\ap\in\Pi_l$
is equal to $2(h^*(\Delta)-2)+2=2h^*(\Delta)-2$.

2. For $\vts$, we have ${\mathcal P}(\vts)=\Delta_s\cup\{0\}$.
Here again ${\frak W}(\vts)$ is a simple graph.
After deleting $\{0\}$ from this weight diagram, we obtain two isomorphic
connected components corresponding to $\Delta^+_s$ and $\Delta^-_s$, etc.
The difference with the previous example is that we loose only two
edges of each {\sl short\/} type after removing $\{0\}$.

3. Let $\g$ be of type $\GR{B}{n}$ and $\lb=\vp_n$.
Here the $\g$-module $\VV_{\vp_n}$ is weight multiplicity free,
$\dim \VV_{\vp_n}=2^n$, and the weights are
$\frac{1}{2}(\pm\esi_1\pm\esi_2 \ldots \pm\esi_n)$. 
The simple root $\ap_n=\esi_n$ can be subtracted from a weight of this
representation if and only if the last coordinate has sign $+$.
Hence the number of edges of type $\ap_n$ in ${\frak W}(\VV_{\vp_n})$
equals $2^{n-1}$. Similarly, one computes that the number of edges of
any other type equals $2^{n-2}$.

4. Let $\g$ be of type $\GR{E}{6}$ and $\lb=\vp_1$. 
Here $\dim\VV_{\vp_1}=27$.
Then ${\mathcal P}(\VV_{\vp_1})$ is a so-called minuscule poset
and ${\frak W}(\VV_{\vp_1})$ has 36 edges. Surely, we 
have 6 edges of each type. 
\end{rem}%
%


\section{The poset of {\sf ad}-nilpotent ideals of $\be$}
\label{pos_all}
\setcounter{equation}{0}

\noindent
There is another natural poset attached to $\Delta^+$, where one can define
the type of edges in the Hasse diagram. As we shall see, this leads to
interesting combinatorial results.
\\[.6ex]
Recall that $\Delta^+$ is equipped with the partial order ``$\curle$''.
A subset $S\subset \Delta^+$ is called an {\it upper ideal\/}, if the conditions
$\gamma\in S$ and $\gamma\curle\tilde\gamma$ imply $\tilde\gamma\in S$.
The geometric counterpart of an upper ideal is an \adn ideal.
A subspace $\ce\subset\be$ is said to be an {\sf ad}-{\it nilpotent ideal\/} (of $\be$), 
if it is contained in $\ut^+$ and satisfies the condition $[\be,\ce]\subset \ce$.
If $\ce$ is an \adn ideal, then $\ce=\underset{\gamma\in I_\ce}{\oplus}\g_\gamma$, 
where $I_\ce$ is a subset of $\Delta^+$. We also say that 
$I_\ce$ is the {\it set of roots\/} of $\ce$.
Obviously, one obtains in this way a bijection between the \adn ideals of $\be$ and
the upper ideals of $\Delta^+$. 

The set of all \adn ideals is denoted by $\AD$.
In view of the above bijection, we may (and will) identify $\AD$ with the set
of all upper ideals of $\Delta^+$.
Whenever we wish to explicitly indicate that $\AD$ depends on $\g$,
we write $\AD(\g)$.
We regard $\AD$ as poset under the usual containment. For instance,
$\ut^+$ (or $\Delta^+$) is the unique maximal element of $\AD$.
It was shown in \cite{cp2} that
\begin{equation}  \label{chislo-mi}
   \#\AD=\prod_{i=1}^n \frac{h+e_i+1}{e_i+1} \ ,
\end{equation}
where $e_1,\ldots,e_n$ are the exponents of $\Delta^+$.
Let $\EAD$ denote the set of edges of the Hasse diagram
$\fH(\AD)$.
Clearly, a pair of \adn ideals $\ce,\ce'$ gives rise to an edge of $\fH(\AD)$
if and only if $\ce'\subset \ce$ and $\dim\ce=\dim\ce'+1$. 
Combinatorially:
$I_{\ce'}\subset I_\ce$ and $\# I_\ce=\# I_{\ce'}+1$. Hence
$I_\ce=I_{\ce'}\cup\{\gamma\}$. It is easily seen that $I_\ce\setminus \{\gamma\}$
is again an upper ideal if and only if $\gamma$ is a generator of
$\ce$ in the sense of the following definition.
\\
An element $\gamma\in I_\ce$ is called
a {\it generator\/} of $\ce$ (or $I_\ce$), if $\gamma-\ap\not\in I$ for any $\ap\in\Pi$.
In other words, $\gamma$ is a minimal element of $I_\ce$ with respect to ``$\curle$''.
We write $\Gamma(\ce)$ for the set of generators of $\ce$.
Hence $\ee=(\ce,\ce')$ is an edge of $\fH(\AD)$ if and only if 
$I_{\ce'}=I_\ce\setminus\{\gamma\}$ for some $\gamma\in \Gamma(\ce)$.
In this case we also say that $\ee$ {\it originates\/} in $\ce'$ and {\it terminates\/} 
in $\ce$  (or, $\ce$ is the {\it terminating ideal\/} 
of $\ee$). Thus, we have proved

\begin{s}{Proposition}  \label{edge_bij}
There is a bijection between $\EAD$ and the set of pairs
$(\ce,\gamma)$, where $\ce\in\Ad$ and  $\gamma\in\Gamma(\ce)$. More precisely, for
any $\ce\in\AD$, the edges terminating in $\ce$ are in a bijection with $\Gamma(\ce)$.
\end{s}%
In \cite[Sect.\,6]{duality}, we introduced, for each simple Lie algebra $\g$,
a generalised Narayana polynomial $\N_{\g}$. By definition,
it is the generating function that counts the \adn ideals with respect to
the number of generators, i.e.,
\[
    \N_{\g}(q)=\sum_{i=0}^{n}\#\{\ce\in \AD(\g) \mid \#\Gamma(\ce)=k\}{\cdot} q^k \ .
\]
In case $\g=\slno$, one obtains the Narayana polynomials
($q$-analogues of the Catalan number).
Obviously, $\N_{\g}(1)=\#\AD(\g)$. 
The following readily follows from the definition of $\N_\g$ and 
Proposition~\ref{edge_bij}.

\begin{s}{Proposition}    \label{nar+edges}
  $ 
\frac{d}{dq}\N_{\g}(q)\vert_{q=1}=\#\EADg$.
\end{s}%
\vskip-1.2ex
\begin{s}{Corollary} \  $\displaystyle
\#\EADg=\frac{n}{2}\#\AD(\g)=\frac{n}{2}\prod_{i=1}^n \frac{h+e_i+1}{e_i+1}$.

\end{s}\begin{proof}
It was observeed in \cite[Sect.\,6]{duality} that $ \N_{\g}(q)=\sum_{i=0}^n a_i q^i
$ is a palindromic polynomial, i.e., $a_j=a_{n-j}$ for all $j$.
But it is easily seen that $f'(1)= \frac{\deg f}{2} f(1)$ for any palindromic
polynomial $f$.
\end{proof}%
Below, we provide the numbers $\#\AD(\g)$ and 
$\#\EADg$ for all simple Lie algebras.
\\[.7ex]
\begin{tabular}{c|cccccccc|}
  & $\GR{A}{n}$ & $\GR{B}{n}$, $\GR{C}{n}$ & $\GR{D}{n}$ & $\GR{E}{6}$ & $\GR{E}{7}$ 
& $\GR{E}{8}$ & $\GR{F}{4}$ & $\GR{G}{2}$ \\ \hline
$\#\AD$   & {\rule{0pt}{3ex}}$\frac{1}{n+2}\genfrac{(}{)}{0pt}{}{2n+2}{n+1}$ &
$\genfrac{(}{)}{0pt}{}{2n}{n}$  & 
$\genfrac{(}{)}{0pt}{}{2n}{n}-\genfrac{(}{)}{0pt}{}{2n-2}{n-1}$  &
833 & 4160 & 25080 & 105 & 8
\\
$\#\EAD$ &  {\rule{0pt}{3.5ex}}$\genfrac{(}{)}{0pt}{}{2n+1}{n+2}$ &
$n\genfrac{(}{)}{0pt}{}{2n-1}{n}$ &
$n(\genfrac{(}{)}{0pt}{}{2n-1}{n}-\genfrac{(}{)}{0pt}{}{2n-3}{n-1})$  &
2499 & 14560 & 100320 & 210 & 8 \raisebox{-1.5ex}{\rule{0pt}{3.5ex}}\\ 
\end{tabular}
\begin{rem}{Remark}  \label{politop}
In \cite{polytope}, a simple convex polytope ${\MM}(\g)$
is associated to an arbitrary irreducible finite root system or simple
Lie algebra (for types $\GR{A}{n}$ and $\GR{B}{n}$ these polytopes were known before).
It is curious that the number of vertices of ${\MM}(\g)$ is
$\#\AD(\g)$, and the number of edges of ${\MM}(\g)$ is
$\#\EADg$. We do not know of whether there is a more deep
connection between ${\MM}(\g)$ and the Hasse diagram 
$\fH(\AD(\g))$. At least, $\gH(\AD(\g))$ is not isomorphic to the graph of
${\MM}(\g)$.
\end{rem}%
In order to define the type of an edge, we need some results on a relationship
between $\AD$ and certain elements of the affine Weyl group.
Let us recall the necessary setup.
\\[.7ex]
We have $V:=\te_{\Bbb R}=\oplus_{i=1}^n{\Bbb R}\ap_i$ and 
$(\ ,\ )$ a $W$-invariant inner product on $V$. As usual,
$\mu^\vee=2\mu/(\mu,\mu)$ is the coroot
for $\mu\in \Delta$. Then $Q^\vee=\oplus _{i=1}^n {\Bbb Z}\ap_i^\vee$  
is the coroot lattice in $V$.
\\[.6ex]
Letting $\widehat V=V\oplus {\Bbb R}\delta\oplus {\Bbb R}\lb$, we extend
the inner product $(\ ,\ )$ on $\widehat V$ so that $(\delta,V)=(\lb,V)=
(\delta,\delta)=
(\lb,\lb)=0$ and $(\delta,\lb)=1$.
Then 
\begin{itemize}
\item[] \ 
$\widehat\Delta=\{\Delta+k\delta \mid k\in {\Bbb Z}\}$ is the set of affine
(real) roots; 
\item[] \ $\HD^+= \Delta^+ \cup \{ \Delta +k\delta \mid k\ge 1\}$ is
the set of positive affine roots; 
\item[] \ $\widehat \Pi=\Pi\cup\{\ap_0\}$ is the corresponding set
of affine simple roots. 
\end{itemize}
Here  $\ap_0=\delta-\theta$.
For $\ap_i$ ($0\le i\le n$), let $s_i$ denote the corresponding 
reflection in $GL(\HV)$.
That is, $s_i(x)=x-2(x,\ap_i)\ap_i^\vee$ for any $x\in \HV$.
The affine Weyl group, $\HW$, is the subgroup of $GL(\HV)$
generated by the reflections $s_i$, $i=0,1,\ldots,n$.
If the index of $\ap\in\widehat\Pi$ is not specified, then we merely write
$s_\ap$. 
The inner product $(\ ,\ )$ on $\widehat V$ is
$\widehat W$-invariant. The notation $\beta>0$ (resp. $\beta <0$)
is a shorthand for $\beta\in\HD^+$ (resp. $\beta\in -\HD^+$).
The length function on $\widehat W$ with respect
to  $s_0,s_1,\dots,s_p$ is denoted by $\ell$.
\\
It was proved by Cellini and Papi that there is a bijection between 
\adn ideals and certain elements of $\HW$, see \cite{cp1}. This can be described as
follows.

Given $\ce\in \AD$ with the corresponding upper ideal $I_\ce\subset\Delta^+$, there is
a unique element $w_{\ce}\in\HW$ satisfying the
following properties:
\begin{enumerate}
\item[($\Diamond$)] \ For $\gamma\in \Delta^+$, we have $\gamma\in I_\ce$ if and only
if $w_\ce(\delta-\gamma) < 0$;
\item[\sf (dom)] \ $w_\ce(\ap)>0$ for all $\ap\in\Pi$;
\item[\sf (min)] \ if   
$w_\ce^{-1}(\ap_i)=k_i\delta+\mu_i$  ($\ap_i\in\HP$), where
$\mu_i\in \Delta$ and $k_i\in \Bbb Z$, then $k_i\ge -1$.
\end{enumerate}
Following Sommers \cite{eric},
the element $w_\ce$ is said to be the {\it minimal element of\/} $\ce$.
The minimal element of $\ce$ can also be characterised as the unique element
of $\HW$ satisfying properties ($\Diamond$), {\sf (dom)}, 
and having the minimal possible length.
This explains the term. 
The elements of $\HW$ satisfying the last two properties are called {\it minimal}.
The set of minimal elements of $\HW$ is denoted by $\HW_{min}$.
One of the main results of \cite{cp1} is that
the correspondence $\ce\mapsto w_\ce$ sets up a bijection between $\AD$ and $\HW_{min}$.
Conversely, if $w\in\HW_{min}$, then $\ce_w$ stands for the corresponding \adn ideal
and $I_w$ is the set of roots of $\ce_w$.
\\[.6ex]
In \cite[Theorem\,2.2]{duality}, a characterisation of the generators of $\ce$ was 
given in terms of $w_\ce$. Namely,
$\gamma\in \Gamma(\ce)$ if and only if $w_\ce(\delta-\gamma)\in -\HP$.
\\
Now, we are ready to define the type of an edge in $\fH(\AD)$.

\begin{rem}{Definition}   \label{def-type}
If $\ee=(\ce,\ce')$ is an edge of $\fH(\AD)$, with
$I_{\ce'}=I_\ce\setminus\{\gamma\}$,
then the {\it type\/} of $\ee$ is the affine simple root $w_\ce(\gamma-\delta)$.
\end{rem}%
Thus, the parameter set for edge types is $\HP$.
Let $\EAD_i$ denote the set of edges of type $\ap_i$, $i=0,1,\ldots,n$.
It is a natural problem to find the number of edges of each type.
We provide two geometric descriptions of $\EAD_i$.
To this end, we recall another bijection due to Cellini and Papi.
Set $D_{min}=\{x\in V \mid (x,\ap)\ge -1 \ \forall\ap\in\Pi \ \ \& \ 
\ (x,\theta)\le 2\}$. It is a simplex in $V$.
\begin{s}{Proposition {\ququ \cite[Prop.\,2 \& 3]{cp2}}}
\label{opis-mi} 
There is a natural bijection between $\AD$ and $D_{min}\cap Q^\vee$.
\end{s}%
In \cite{cp2}, this bijection was established using the isomorphism
$\HW\simeq W\ltimes Q^\vee$ and the affine-linear action of $\HW$ on $V$.
This can also be explained entirely in terms of the linear action of
$\HW$ on $\HV$. 
If $w_\ce^{-1}(\ap_i)=\mu_i+k_i\delta, \quad i=0,1,\ldots,n$, then
we define the point $z_\ce\in V$ by the equalities
$(\ap_i, z_\ce)=k_i$, $i=1,\ldots,n$. (Since $\ap_0=\delta-\theta$
and $\delta$ is $\HW$-invariant, this implies
$(\theta,z_\ce)=1-k_0$.)
Since $w_\ce\in\HW_{min}$, we have $k_i\ge -1$, hence $z_\ce\in D_{min}$.
The correspondence  $\ce\mapsto z_\ce$ is the required bijection.

Set $F_i=\{x \in D_{min} \mid (x,\ap_i)=-1\}$, $i\ge 1$,
and $F_0=\{x \in D_{min}\mid (x,\theta)=2\}$.
These are all the facets of $D_{min}$.

\begin{s}{Theorem}  \label{tip-edges}
We have $\#\EAD_i=\#(F_i\cap Q^\vee)$, i.e.,
the number of edges in $\EAD_i$ is equal to the number of points $z_\ce$
lying in $F_i$ ($i=0,1,\ldots,n$).
\end{s}\begin{proof}
This is a simple combination of preceding results.
Clearly, in place of edges of type $\ap_i$ one can describe their terminating 
ideals.
Let $\ee=(\ce,\ce')$ be an edge of type $\ap_i$, where $I_{\ce'}=I_\ce\setminus
\{\gamma\}$. Then 
$w_\ce(\delta-\gamma)=-\ap_i$. Hence $w_\ce^{-1}(\ap_i)=\gamma-\delta$.
Comparing with the above definition of $z_\ce$ shows that $z_\ce\in F_i$.
Conversely, if $z_\ce$ lies in $F_i$, then 
$w_\ce(\delta-\gamma)=\ap_i$ for some $\gamma\in I_\ce$.
By \cite[Theorem\,2.2]{duality}, such $\gamma$ is necessarily a generator of $\ce$,
and $(I_\ce,I_\ce\setminus\{\gamma\})$
gives rise to an edge of type $\ap_i$.
\\
Thus, we have proved that there is an edge of type $\ap_i$ terminating in $\ce$
if and only if $z_\ce\in F_i$.
\end{proof}%
Yet another characterisation of the type of an edge can be given using alcoves.
Recall that 

the (open) {\it dominant Weyl chamber\/} is
${\gC}=\{x\in V\mid (x,\ap)>0 \ \ \forall \ap\in\Pi\}$ \ and 

the {\it fundamental alcove\/} is
${\gA}=\{x\in V\mid (x,\ap)>0 \ \ \forall \ap\in\Pi \ \ \& \ (x,\theta)<1\}$.
\\[.6ex]
For $\gamma\in\Delta^+$ and $k\in\Bbb Z$, we set 
$\gH_{\gamma,k}=\{x\in V\mid (\gamma,x)=k\}$. It is an affine hyperplane in $V$.
The connected components of $V\setminus \cup_{\gamma,k}\gH_{\gamma,k}$ 
are called {\it alcoves\/}. 
As is well-known, $\gA$ is one of them, and all alcoves are congruent to $\gA$,
see e.g. \cite{hump}.
The walls of $\gA$ are $\gH_{\ap_i,0}$, $\ap_i\in\Pi$, and $\gH_{\theta,1}$.
Each wall of $\gA$ is equipped with the "type", which is an element of $\HP$.
Namely, $\gH_{\theta,1}$
is the wall of type $\ap_0$ and $\gH_{\ap,0}$ ($\ap\in\Pi$) is the wall
of type $\ap$. 
\\[.6ex]
Given $\ce\in\AD$, we wish to determine the types of edges 
terminating in $\ce$. To this end, 
consider the alcove $w_\ce^{-1}\ast\gA$. Here `$\ast$' stands for the 
affine-linear action of $\HW$ on $V$. 
Since $w_\ce$ satisfies condition {\sf (dom)}, we have 
$w_\ce^{-1}{\ast}\gA \subset \gC$, see~\cite{cp1}. 

\begin{s}{Theorem}   \label{geom}
In the above setting, there is a bijection between the edges of $\EAD$ terminating
in $\ce$ and the walls of $w_\ce^{-1}{\ast}\gA$
separating this alcove from the origin. 
If $H$ is such a wall of $w_\ce^{-1}{\ast}\gA$, then the type of the corresponding
edge coincides with the type of the wall $w_\ce{\ast} H$ of $\gA$.
\end{s}\begin{proof}
{\it First}, suppose $H=\gH_{\gamma,k}$ separates $w_\ce^{-1}{\ast}\gA$ from the 
origin. Then $w_\ce(k\delta-\gamma)<0$ \cite[Eq.\,(1.1)]{cp1}. 
Let us show that $H$ gives rise to an edge terminating in $\ce$.
If $w_\ce{\ast}\gH_{\gamma,k}=\gH_{\ap_i,0}$, then  using \cite[Eq.\,(3.3)]{norm-p}
we obtain $w_\ce(k\delta-\gamma)=-\ap_i$. Since $w_\ce$ is minimal,
the property {\sf (min)} forces $k=1$. Hence $\gamma\in \Gamma(\ce)$
and the respective edge is of type $\ap_i$.
If $w_\ce{\ast}\gH_{\gamma,k}=\gH_{\theta,1}$, then we obtain in a similar way
that $k=1$ and $w_\ce(k\delta-\gamma)=-\ap_0$. Hence $H$ gives rise to an
edge of type $\ap_0$.

{\it Conversely\/}, suppose $\gamma\in \Gamma(\ce)$ corresponds to an edge of type
$\ap\in\HP$. Then arguing backwards, we show that $H=\gH_{\gamma,1}$
separates $w_\ce^{-1}{\ast}\gA$ from the 
origin and $w_\ce{\ast}H$ is the wall of $\gA$ of type $\ap$.
\end{proof}%
Finally, we show how to determine the edges terminating in $\ce$ in terms of 
reduced decompositions of $w_\ce$.
Suppose that $\ell(w_\ce)=m$, and $w_\ce=s_{i_1}\ldots s_{i_m}$ is a reduced 
decomposition. It follows from property {\sf (dom)} that $\ap_0$ is the only
simple root that is made negative by $w_\ce$. Hence $s_{i_m}=s_0$.
But at the other side we may have several possibilities.

\begin{s}{Proposition}   \label{reduced}
Given $\ce$, there is  an edge of type $\ap_i$ terminating in $\ce$
if and only if there is a reduced decomposition of $w_\ce$ starting with
$s_i$.
\end{s}\begin{proof}
Indeed, if $w_\ce=s_i w'$ and $\ell(w')=\ell(w_\ce)-1$, then
$w_\ce^{-1}(\ap_i)<0$. Using properties {\sf (min)} and {\sf (dom)}, we see
that this implies that $w_\ce^{-1}(\ap_i)=\gamma-\delta$ for some 
$\gamma\in\Delta^+$. Hence $\gamma\in\Gamma(\ce)$ and the corresponding
edge is of type $\ap_i$. This argument can be reversed.
\end{proof}%
It is not clear as yet how to uniformly compute the number of edges of each type.
At least, these numbers can be different for simple roots of the same length.

\begin{rem}{Example}   \label{ex3.1}
For $\g={\frak sl}_4$, the Hasse diagram of $\AD$ has 21 edges. Their distribution 
with respect to type is $(5,\,6,\,4,\,6)$, see Figure~\ref{hasse:sl4}.
The asterisks point out the generators of ideals.
\begin{figure}[htb]
\begin{center}
\setlength{\unitlength}{0.02in}
\begin{picture}(300,113)(-5,0)
\Yboxdim2.5mm
\put(0,50){$\varnothing$}
\put(30,50){\young(\ast )}
\put(60,75){\young(\ast \ )}
\put(58,25){\young(:\ ,:\ast)}
\put(90,100){\young(\ast \ \ )}
\put(90,47){\young(:\ast \ ,::\ast)}
\put(88,0){\young(::\ ,::\ ,::\ast )}
\put(130,97){\young(\ast \ \ ,::\ast )}
\put(130,47){\young(:\ \ ,:\ast \ )}
\put(130,0){\young(:\ast\ ,::\ ,::\ast )}
\put(170,97){\young(\ast \ \ ,:\ast \ )}
\put(170,45){\young(\ast \ \ ,::\ ,::\ast )}
\put(170,0){\young(:\ \ ,:\ast \ ,::\ast )}
\put(215,45){\young(\ast \ \ ,:\ast \ ,::\ast )}
\put(231,50){=$\ut^+$}
  \put(7,52){\vector(1,0){20}} \put(14,54){\tiny $\ap_0$}
\put(38,54){\vector(1,1){20}} \put(45,68){\tiny $\ap_3$}
\put(38,50){\vector(1,-1){20}} \put(47,41){\tiny $\ap_1$} 
\put(71,80){\vector(1,1){18}} \put(77,93){\tiny $\ap_2$}  %
\put(71,75){\vector(1,-1){20}} \put(79,68){\tiny $\ap_1$} 
\put(71,30){\vector(1,1){20}} \put(75,41){\tiny $\ap_3$}  
\put(71,25){\vector(1,-1){20}} \put(77,21){\tiny $\ap_2$}  %
\put(107,102){\vector(1,0){20}} \put(114,105){\tiny $\ap_2$}  
\put(107,52){\vector(1,0){24}} \put(118,54){\tiny $\ap_0$}
\put(107,5){\vector(1,0){23}} \put(114,0){\tiny $\ap_2$}  
\put(107,56){\vector(3,4){27}} \put(117,78){\tiny $\ap_1$} 
\put(107,47){\vector(3,-4){25}} \put(118,36){\tiny $\ap_3$}  
\put(147,102){\vector(1,0){21}} \put(154,105){\tiny $\ap_1$} 
\put(147,5){\vector(1,0){24}} \put(154,0){\tiny $\ap_3$}  
\put(147,56){\vector(3,4){27}} \put(152,59){\tiny $\ap_0$} 
\put(147,47){\vector(3,-4){25}} \put(152,41){\tiny $\ap_0$}
\put(147,98){\vector(3,-4){25}} \put(154,91){\tiny $\ap_3$}  
\put(147,9){\vector(3,4){27}}  \put(151,23){\tiny $\ap_1$} 
\put(187,98){\vector(3,-4){25}} \put(196,90){\tiny $\ap_1$} 
\put(188,52){\vector(1,0){23}}  \put(196,54){\tiny $\ap_0$}
\put(187,9){\vector(3,4){27}} \put(196,16){\tiny $\ap_3$} 
\end{picture}
\caption{The Hasse diagram of $\AD({\frak sl}_4)$ with types of the edges} \label{hasse:sl4}
\end{center}
\end{figure}
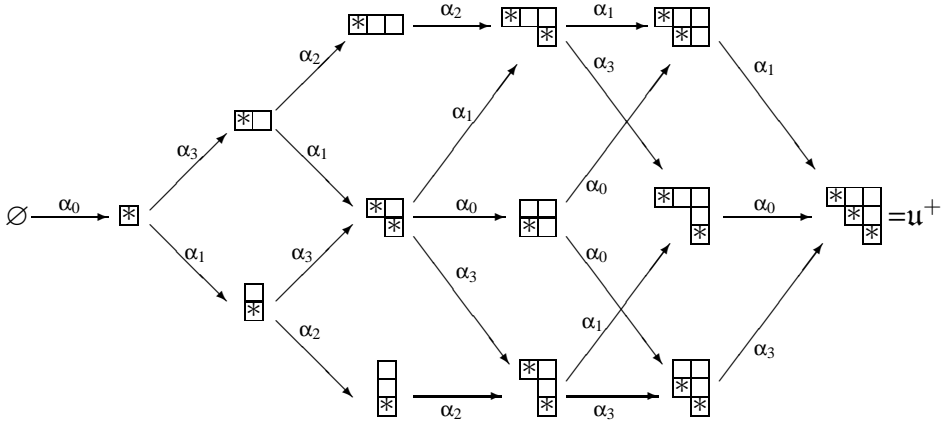
\end{rem}%
To compute the type of edges, we use explicit expressions for the minimal
elements of all \adn ideals. For instance, if $\ce=\ut^+$, then
the corresponding minimal element is 
$s_0s_1s_3s_0s_3s_1s_2s_1s_3s_0$. It is not hard to find reduced expressions of it starting
with $s_1$ or $s_3$.


\section{The poset of Abelian ideals of $\be$}
\label{pos_ab}
\setcounter{equation}{0}

\noindent
A subspace $\ce\subset\be$ is called an {\it Abelian ideal\/}, if 
$[\be,\ce]\subset \ce$ and $[\ce,\ce]=0$. (It is easily seen that any
Abelian ideal is \adn, i.e., it is contained in $\ut^+$.)
In terms of roots, $\ce$ is Abelian if and only if 
$(I_\ce+I_\ce)\cap \Delta^+=\varnothing$.
The respective upper ideal $I_\ce$ is said to be Abelian, too.
The sub-poset of $\AD$ consisting of all Abelian ideals is denoted by $\Ab$
or $\Ab(\g)$. Below, the symbol $\ah$ is used to denote an Abelian ideal.

We keep the notation of the previous section.
In particular, to any $\ah\in\Ab$ we attach the set of roots $I_\ah$, the
minimal element $w_\ah$, and the alcove $w_\ah^{-1}\ast\gA \subset \gC$.
Following Peterson, an element
$w\in \HW_{min}$ is said to be {\it minuscule\/}, if the corresponding
\adn ideal is Abelian. Write $\HW_{mc}$ for the set of all
minuscule elements. 
We are going to consider the edges of the Hasse diagram
$\fH(\Ab)$ and their types. In this case, our understanding
of the situation is better than that in Section~\ref{pos_all}.

A nice result of D.\,Peterson asserts that $\#\Ab(\g)=2^{\rk\g}$
for any simple Lie algebra $\g$. (See \cite{ko1} and \cite{cp1} for various
proofs). The proof of Cellini and Papi is based on the observation
that $\ce\in\AD$ is Abelian if and only if $w_\ce^{-1}\ast\gA
\subset 2\gA$. It turns out that the number of edges depends only on the rank, too.

\begin{s}{Theorem}   \label{main:ab}
If\/ $\rk\g=n$, then \ $\#\EAbg= (n+1)2^{n-2}$.
\end{s}\begin{proof}
Recall that the hyperplanes $\gH_{\gamma,k}$ with $\gamma\in\Delta^+$
and $k\in \Bbb Z$ cut $V$ into (open) alcoves congruent to $\gA$.
In particular, the `big' simplex $2\gA$ contains $2^n$ alcoves.
By \cite{cp1}, the alcoves in $2\gA$ bijectively correspond to the Abelian ideals,
via the mapping $(\ah\in \Ab) \mapsto (w_\ah^{-1}\ast\gA \subset 2\gA)$.
Hence the number of Abelian ideals is $2^n$.
Consider two alcoves inside $2\gA$ that have a common wall.
By Theorem~\ref{geom}, this wall gives rise to an edge of $\gH(\Ab)$.
Hence the number of edges is equal to the total number of internal walls
between alcoves inside $2\gA$.
It is an easy exercise to compute this number.
(The total number of walls of all these alcoves is $(n+1)2^n$;
the number of alcove walls on the boundary of $2\gA$ is $(n+1)2^{n-1}$;
all the remaining walls are counted twice. Hence the answer.)
\end{proof}%
The hyperplanes $\gH_{\gamma,k}$
cut $2\gA$ into $2^n$ alcoves. It is natural to ask which hyperplanes do meet
$2\gA$? A partial answer given in \cite{cp1} says that if 
$2\gA\cap \gH_{\gamma,k}\ne\varnothing$, then $k=1$. Therefore, it remains
to characterise the possible roots $\gamma$. This leads to the following
definition.

\begin{rem}{Definition}   \label{def:ab_root}
Let $\gamma\in\Delta^+$. We say that $\gamma$ is {\it commutative\/} if
the upper ideal generated by $\gamma$ is Abelian. This Abelian ideal is denoted
by $\ah(\gamma)$.
\end{rem}%
It is easily seen that $\gamma$ is commutative if and and only if
there is an Abelian ideal $\ah$ such that $\gamma\in\Gamma(\ah)$.
Clearly, the set of all commutative roots forms an \adn ideal.
This \adn ideal is Abelian if and only if there is a unique maximal Abelian
ideal. (Although this is not used in the sequel, we note that this
happens only for $\GR{C}{n}$ and $\GR{G}{2}$.)

\begin{s}{Lemma}   \label{cuts}
$\gH_{\gamma,1}$ meets $2\gA$ if and only if $\gamma$ is commutative.
\end{s}\begin{proof}
1. Suppose that $\gamma$ is not commutative. Then there exist 
$\nu_1,\nu_2\curge\gamma$ such that $\nu_1+\nu_2=\theta$.
If $x\in\gC\cap \gH_{\gamma,1}$, then $(x,\nu_i)\ge 1$.
Hence $(x,\theta)\ge 2$, i.e., $x\not\in 2\gA$.

2. Suppose $\gamma$ is commutative, and let $w\in\HW_{mc}$ be the minimal element
of $\ah(\gamma)$. Then $w^{-1}\ast\gA\subset 2\gA$ and $H_{\gamma,1}$ is the wall
separating $w^{-1}\ast\gA$ from the origin.
\end{proof}%
Next, we list some properties of (non)-commutative roots in
$\Delta^+$. It would be very interesting to find an a priori proof.

\begin{s}{Theorem}  \label{number:com} 
\begin{itemize}
\item[\sf (i)] \ If $\theta=\sum_{i=1}^nc_i\ap_i$, then
$[\theta/2]:=\sum_{i=1}^n [c_i/2]\ap_i$ is the unique maximal non-commutative
root. (In case of $\slno$, we temporarily assume that $0$ is a root.)
\item[\sf (ii)] \
a) \ If the Dynkin diagram of $\Delta$ has no branching nodes, then
the number of commutative roots is \ $n(n+1)/2$. 
\\
b) \ If there is a branching node and $n_1,n_2,n_3$ are the lengths
of tails of the diagram obtained after deleting the branching node, 
hence $n_1+n_2+n_3=n-1$, then 
the number of commutative roots is \ $n(n+1)/2+ n_1n_2n_3$.
\end{itemize}
\end{s}\begin{proof}
(i) \ Suppose $\gamma=\sum_{i=1}^n m_i\ap_i\in\Delta^+$ and $m_i> [c_i/2]$
for some $i$. Then $\gamma$ is obviously commutative. Hence all non-commutative
roots satisfy the constraint $m_i\le [c_i/2]$ for all $i$. 
Then one observes a mysterious fact that both $[\theta/2]$ and 
$\theta-[\theta/2]$ are
always roots. Since $[\theta/2]\curle\theta-[\theta/2]$,
we see that $[\theta/2]$ is non-commutative.

(ii) \ 
1. For $\g=\slno$, all positive roots are commutative. \\
2. For $\g=\spn$,  the commutative roots are $\{\esi_i+\esi_j\mid 1\le i\le j \le n\}$.\\
3. For $\g=\sono$,  the commutative roots are 

$\{\esi_i+\esi_j\mid 1\le i<j \le n\}\cup
\{\esi_1-\esi_i\mid 2\le i\le n\}\cup\{\esi_1\}$. \\
4. For $\g=\sone$, the commutative roots are 

$\{\esi_i+\esi_j\mid 1\le i<j \le n\}\cup
\{\esi_1-\esi_i\mid 2\le i\le n\}\cup\{\esi_i-\esi_n\mid 2\le i\le n-1\}$. \\
5. For the exceptional Lie algebras, one can perform explicit verification.
E.g., the number of commutative roots is equal to
$25,34,44$ for $\GR{E}{6},\GR{E}{7},\GR{E}{8}$, respectively. 
\end{proof}%
{\bf Remarks.} {\sf 1.} For each $\gamma\in\Delta^+$, it is true that
$[\gamma/2]\in \Delta^+\cup\{0\}$ and $\gamma-[\gamma/2]\in\Delta^+$.

{\sf 2}. Comparing Lemma~\ref{cuts} and Theorem~\ref{number:com}(ii) shows 
that the number of hyperplanes that is needed 
to cut $2\gA$ into $2^n$ congruent simplices depends not only on $n$, but also 
on the angles between walls of $\gA$. 
\\[.7ex]
In the Abelian case, the relationship between minimal elements
and \adn ideals works much better. This is explained by the fact
that for $\ce\in\AD$ we have 
$\ell(w_\ce)= \dim\ce$ if and only if $\ce\in\Ab$.
For, it is known in general that
$\ell(w_\ce)=\dim\ce+\dim[\ce,\ce]+\dim [\ce,[\ce,\ce]]+\ldots$, see
\cite{cp1}.

Suppose that $w\in \HW_{mc}$ and $w^{-1}(\ap)<0$ for $\ap\in\HP$.
Then the element $w'=s_\ap w$ whose length is one less
is again minuscule. 
Conversely, if $w'\in\HW_{mc}$ and $w'(\delta-\gamma)=\ap\in\HP$ for some
$\gamma\in\Delta^+$, then $w=s_\ap w'$ is also minuscule and 
$I_{w}=I_{w'}\cup\{\gamma\}$ is the set of roots of an Abelian
ideal \cite[Theorem\,2.4]{lp}.
The procedure of such {\it elementary extensions\/} of Abelian ideals was studied in
\cite{lp}. Properties of elementary extensions
have useful consequences for describing the types of edges
in $\fH(\Ab)$.
The situation here resembles that in Section~\ref{odin}, i.e.,
the number of edges of any given type depends only on the length
of the respective affine simple root. 

Let us look again at Definition~\ref{def-type} (of the type of an edge).
If $\gamma\in\Gamma(\ce)$ and $I_{\ce'}:=I_\ce\setminus\{\gamma\}$, then
$w_\ce(\gamma-\delta)=:\ap$ is the type of $(\ce,\ce')$.
In this situation, we say that the generator $\gamma$ is {\it of class\/} $\ap\in\HP$
in the ideal $\ce$. This notion is not, however, very convenient in general, since
$\gamma$ can have another class in another \adn ideal in which it is a generator.
This is already seen in case $\g={\frak sl}_4$,
see Figure~\ref{hasse:sl4}.
But in the Abelian case this unpleasant phenomenon does not occur.

\begin{s}{Theorem}   \label{class:ab}
Let $\gamma$ be a commutative root. Then the class of $\gamma$ does not depend on
an Abelian ideal in which $\gamma$ is a generator.
\end{s}\begin{proof}
Recall that $\ah(\gamma)$ is the minimal (Abelian)
\adn ideal generated by $\gamma$.
All other Abelian ideals are obtained from $\ah(\gamma)$ via sequences of elementary 
extensions. So, it is enough to prove that if $\ah'\subset\ah$ are Abelian, 
$\dim\ah'=\dim\ah-1$, and $\gamma\in\Gamma(\ah)\cap\Gamma(\ah')$, then
the classes of $\gamma$ in $\ah$ and $\ah'$ are equal.

We have $I_{\ah}=I_{\ah'}\cup\{\mu\}$ for some $\mu$. Since $\mu$ is a generator
of $\ah$, we have $w_\ah(\mu-\delta)=\ap_i\in\HP$ and $w_\ah=s_iw_{\ah'}$.
Similarly, $w_\ah(\gamma-\delta)=\ap_j\in\HP$. 
As $\gamma$ remains a generator of $\ah'$, we have $\gamma\not\curle\mu$
and $\mu\not\curle\gamma$. In particular, $(\gamma,\mu)=0$.
Hence $(\ap_i,\ap_j)=(\mu-\delta,\gamma-\delta)=(\mu,\gamma)=0$.
Therefore 

$w_{\ah'}(\gamma-\delta)=s_iw_\ah(\gamma-\delta)=s_i(\ap_j)=
\ap_j=w_{\ah}(\gamma-\delta)$.
\end{proof}%
Thus, for any commutative root, the notion of class is well-defined.
We write $cl(\gamma)$ for the class of $\gamma$. 
The class can be regarded as the map \\
$cl: \text{\{commutative roots\}} \to \HP$.
One can use the  
following strategy for computing the cardinality of $\EAb_i$:

{\bf --} \ Determine  the class of each commutative root;

{\bf --} \ If $cl^{-1}(\ap_i)=:\{\mu_1,\ldots,\mu_t\}$, then for each $\mu_j$ 
one counts the number of Abelian ideals having $\mu_j$ as a generator.
The total number of such ideals equals $\#\EAb_i$.

The following lemma is helpful for explicit computations of the class. 

\begin{s}{Lemma}   \label{sosedi}
\\
1. Let $\gamma,\gamma'$ be commutative roots that 
are adjacent in $\fH(\Delta^+)$. Then $cl(\gamma),\,cl(\gamma')$ are adjacent roots
in the extended Dynkin diagram.
\\
2. Let $\gamma$ be a maximal root in the set $\{\mu\in\Delta^+\mid (\mu,\theta)=0\}$.
Then $cl(\gamma)=\ap_0$.
\\
3. If\/ $\gamma_1,\gamma_2\in \Gamma(\ah)$ for some $\ah\in\Ab$, 
then $(cl(\gamma_1),cl(\gamma_2))=0$.
\end{s}\begin{proof}
1. Assume that $\gamma\curle\gamma'$. 
Let $\ah=\ah(\gamma)$ be the Abelian ideals generated by $\gamma$. Then 
$\gamma'\in \Gamma(I_\ah\setminus\{\gamma\})$.
Suppose that $w_\ah(\gamma-\delta)=\ap\in\HP$, i.e., $cl(\gamma)=\ap$.
As follows from the above discussion on elementary extensions of 
Abelian ideals, $w'=s_\ap w_\ah$ is the minimal element of the upper ideal 
$I_\ah\setminus\{\gamma\}$.
Hence, by Theorem~\ref{class:ab}, $cl(\gamma')=w'(\gamma'-\delta)=:\ap'$.
That is,  $w_\ah(\gamma'-\delta)=s_\ap(\ap')$.
Since $\gamma'\not\in\Gamma(\ah)$, we obtain $s_\ap(\ap')\not\in\HP$.
Hence $(\ap,\ap')\ne 0$.

2. Let $I$ be the set of roots of $\ah(\gamma)$, and set $I'=I\setminus\{\gamma\}$.
Then $I'$ is the set of roots of an Abelian ideal $\ah'$ and
$(\mu,\theta)>0$ for each $\mu\in I'$. By Theorem~4.3 in \cite{lp}, 
we then have $w_{\ah'}=ws_0$ for some $w\in W$. 
Since $I'\to I$ is an elementary extension of Abelian ideals,
we must have $w_{\ah'}(\delta-\gamma)\in \HP$. By the assumption,
$w_{\ah'}(\delta-\gamma)=\delta-w(\gamma)$. Hence the only possibility for
this to be a simple root is $\ap_0$. Thus
$w_\ah=s_0 w_{\ah'}$.

3. Suppose $w_\ah=s_1w'=s_2w''$, where $\ell(w')=\ell(w'')=\ell(w)-1$,
and $w_\ah(\delta-\gamma_i)=-\ap_i$, $i=1,2$. If $(\ap_1,\ap_2)\ne 0$, then
$\ap_1+\ap_2\in\HD^+$. Hence $w_\ah^{-1}(\ap_1+\ap_2)=\gamma_1+\gamma_2-2\delta$.
This means that $\gamma_1+\gamma_2$ is a root, i.e., $\ah$ is not Abelian.
\end{proof}%
%
%
By Theorem~\ref{main:ab}, $\#\EAb$ is always divisible 
by $n+1$. This suggests a tempting idea that it might be true that 
the number of edges of each type equals $2^{n-2}$. But, one immediately finds
that for $\spn$ this is not the case. However, this is the only exception.
Our results on the numbers $\#\EAb_i$ are obtained via case-by-case considerations.
To this end, one has to know the class of each commutative root.
For $\slno$, such an information is essentially presented in \cite[Subsection\,6.1]{lp}
under the name of "filling the Ferrers diagram".
Similar information for the other classical series 
is easily obtained via direct computations, using Proposition~\ref{reduced}
and Lemma~\ref{sosedi}.

\begin{s}{Theorem}  \label{edges:sl}
If $\g=\slno$, $n\ge 2$, then $\#\EAb_i=2^{n-2}$ for each $i\in\{0,1,\ldots,n\}$.
\end{s}{\sl Proof.}\quad
In $\slno$, the positive (=\,commutative) roots are 
$\nu_{ij}:=\esi_i-\esi_j$, $1\le i <j \le n+1$.
The numbering of simple roots is such that $\ap_i=\nu_{i,i+1}$, $i=1,\ldots,n$.
Then $cl(\nu_{ij})=\ap_k$, where $i+j-1\equiv k \pmod{n+1}$.
We work with the usual matrix realisation of $\be$ as the set of 
upper-triangular matrices.
The set of positive roots is identified with the right-justified Ferrers diagram
with row lengths $(n,n-1,\ldots,1)$, and the Abelian upper
ideals are identified with the right-justified subdiagrams of it that fit inside the 
rectangle of shape $j\times (n+1-j)$ for some $j$.
Such subdiagrams are said to be Abelian Ferrers diagrams.
The class of roots is constant along the diagonals parallel to the
antidiagonal (in the $n+1$ by $n+1$ matrix).
In particular, the roots of class $\ap_0$ are exactly the roots on the antidiagonal,
see Figure~\ref{ab:sl}.

\begin{figure}[htb]
\begin{center}
\setlength{\unitlength}{0.02in}
\Yboxdim3ex
{\young(\ \apxx\apxxx\apxxxx\apxxxxx\apo,\ \ \apxxxx\apxxxxx\apo\apx,\ \ \ \apo\apx\apxx,\ \ \ \ \apxx\apxxx,\ \ \ \ \ \apxxxx,\ \ \ \ \ \ )}
\caption{Classes of the commutative roots for ${\frak sl}_6$} \label{ab:sl}
\end{center}
\end{figure}
\noindent
The automorphism group of $\fH(\Ab(\slno))$ (as undirected graph!) is equal to 
$D_{n+1}$, the dihedral group of order $2(n+1)$, see \cite{suter_ejc}. 
Furthermore, Suter explicitly describes a certain automorphism 
$\tau$ of order $n+1$, via a ``sliding procedure in south-east direction''. 
Comparing that procedure with the formula for
the class of roots, and thereby the types of edges in $\EAb$,
one finds that $\tau$ always takes an edge of type $\ap_i$ in $\EAb$ to
an edge of type $\ap_{i+2}$ (with the cyclic ordering of affine simple roots),
see Lemma~\ref{tau} below.
Hence, if $n+1$  odd, then $\tau$ acts transitively on the set
of types of edges,  which completes the proof in this case.
If $n+1$ is even, then $\tau$ has 2 orbits on the set of types,
i.e., $\{0,2,4,\ldots,n-1\}$ and $\{1,3,\ldots,n\}$.
Therefore it is enough to check that the number of Abelian ideals having a
generator of class $\ap_0$ equals $2^{n-2}$. Note that the verification 
performed below does not exploit the hypothesis that $n+1$ is even.

Take the commutative root of class $\ap_0$ that lies in the $k$-th row of the 
Ferrers diagram (matrix), i.e., $\nu_{k,n+2-k}$.
Then $k\le (n+1)/2$.
Consider the set of Abelian Ferrers diagram having this root as
a generator (=\,south-west corner) and the first row of length
$x+k$. 
Using the graphical presentation of ideals, see Figure~\ref{graphic}, we obtain that
the cardinality of this set equals
\\[.8ex]  \hspace*{4ex}
$\left\{\parbox{183pt}{the number of all Ferrers diagrams that fit in 
the rectangle of shape $x\times (k-2)$}\right\}\ \times \   
\left\{\parbox{183pt}{the number of all Ferrers diagrams that fit in 
the rectangle of shape $(n+1-2k-x)\times (k-1)$}\right\}$. 
\begin{figure}[htb]
\setlength{\unitlength}{0.02in}
\begin{center}
\begin{picture}(120,125)(-3,0)
\multiput(0,0)(120,0){2}{\line(0,1){120}}
\multiput(0,0)(0,120){2}{\line(1,0){120}}

\qbezier[80](120,0),(60,60),(0,120)  
\qbezier[80](0,0),(60,60),(120,120)  
\put(85,85){\line(1,0){5}}
\put(85,85){\line(0,1){5}}
\put(85,90){\line(1,0){5}}
\put(90,85){\line(0,1){5}}
\put(65,115){\line(1,0){5}}
\put(65,115){\line(0,1){5}}
\put(65,120){\line(1,0){5}}
\put(70,115){\line(0,1){5}}
%
%
\qbezier[40](65,55),(92,55),(120,55)
\qbezier[40](65,55),(65,85),(65,115)
%
%
\qbezier[15](90,85),(90,70),(90,55)
\qbezier[25](90,85),(105,85),(120,85)
\put(71,100){{\small $(\ast)$}}
%
%
\qbezier[12](90,90),(90,105),(90,120)
\qbezier[12](65,90),(75,90),(85,90)
\qbezier[12](85,90),(85,102),(85,115)
\qbezier[10](70,115),(77,115),(85,115)
\put(100,70){{\small $(\ast)$}}

\put(65,120){$\overbrace%
{\mbox{\hspace{20\unitlength}}}^{x}$}
\put(90,120){$\overbrace%
{\mbox{\hspace{30\unitlength}}}^{k-1}$}

\put(128,100){{\footnotesize $k$}}
\put(128,70){{\footnotesize $n{+}1{-}2k{-}x$}}

\end{picture}
\lefteqn{\raisebox{100\unitlength}%
{$\left. {\parbox{1pt}{\vspace{34\unitlength}}}
\right\}$}}%
\raisebox{68\unitlength}%
{$\left. {\parbox{1pt}{\vspace{29\unitlength}}}
\right\}$}

\end{center}
\caption{
}  
\label{graphic}
\end{figure}

\noindent
(Two rectangles in Figure~\ref{graphic} that have to be filled with Ferrers diagrams 
are marked with $(\ast)$.)
Here the first (resp. second) number equals 
$\genfrac{(}{)}{0pt}{}{x+k-2}{k-2}$ (resp. $\genfrac{(}{)}{0pt}{}{n-k-x}{k-1}$).
Therefore the total number of 
Abelian ideals having
the generator of class $\ap_0$ lying in the $k$-th row is equal to
\[
   \sum_{x=0}^{n+1-2k}\genfrac{(}{)}{0pt}{}{x+k-2}{k-2}\genfrac{(}{)}{0pt}{}{n-k-x}{k-1}
\ .
\]
Using Lemma~\ref{need} below, we obtain this sum equals 
$\genfrac{(}{)}{0pt}{}{n-1}{2k-2}$. Hence the total number of 
Abelian ideals having
a generator of class $\ap_0$ is equal to \\
\hbox to \textwidth{\phantom{$x$} \hfil
  $\displaystyle\sum_{k=1}^{(n+1)/2}\genfrac{(}{)}{0pt}{}{n-1}{2k-2}=2^{n-2}$.
\hfil $\Box$} 
\\[.8ex]
Now we prove two results that have been used in the previous proof.

\begin{s}{Lemma}  \label{tau}
Let $\tau$ be Suter's automorphism of the undirected graph $\fH(\Ab(\slno))$
(to be defined below). 
Then $\tau$ takes an edge of type $\ap_i$ to an edge of type $\ap_{i+2}$, with the cyclic
ordering of the affine simple roots.
\end{s}\begin{proof*}
In this proof, we write $(a,b)$ for the root $\nu_{ab}$.
Let $\ah$ be the Abelian ideal with generators
$\Gamma(\ah)=\{(a_1,b_1),\ldots,(a_k,b_k)\}$. Then 
$1\le a_1<\ldots <a_k<b_1 <\ldots <b_k\le n+1$.
The generators of $\tau(\ah)$ are defined by the following rule (it is a formal version
of the diagrams depicted in \cite{suter_ejc}):

$\bullet$ \ If $b_k=n+1$, then $(a_k,b_k)$ disappears. In all other cases,
$(a_j,b_j)$ is replaced with $(a_j+1,b_j+1)$;

$\bullet$ \ If $a_k+1< b_1$, then the new generator $(1,a_k+2)$ emerges.
\\[1ex]
Now, we have to keep track of all edges incident to $\ah$, not only those terminating
in $\ah$. The edges terminating in $\ah$ correspond to the generators, i.e.,
the south-west corners of the respective Ferrers diagram,
whereas the edges originating in $\ah$ correspond to the maximal roots in
$\Delta^+\setminus I_\ah$
(modulo the constraint that $\ah$ plus the respective root space still yields 
an Abelian ideal). Hence, the full collection of edges incident to $\ah$ is 
determined by the roots: 
\\[.7ex]
terminating in $\ah$: $(a_1,b_1),\ldots,(a_k,b_k)$; 
\\[.7ex]
\phantom{t}originating in $\ah$: \parbox[t]{320pt}{
$(a_1+1,b_2-1),\ldots,(a_{k-1}+1,b_k-1)$; furthermore,
if $a_k+1<b_1$, then the roots $(1,b_1-1),(a_{k}+1,n+1)$ are also needed.}
\\[.6ex]
It is easy to see how these roots transform under $\tau$. If a root $(a,b)$
does not belong to the last column of the matrix, then it merely
goes to $(a+1,b+1)$. Hence the "increment" in the type number is 2.
All possibilities for the roots lying in the last column are easily
handled in a case-by-case fashion.
\end{proof*}%
\begin{s}{Lemma}  \label{need} \\
For $a,b,c\in \Bbb N$ and $a\ge b+c$, we have $\displaystyle
\sum_{x=c}^{a-b}\genfrac{(}{)}{0pt}{}{x}{c}\genfrac{(}{)}{0pt}{}{a-x}{b}=
\genfrac{(}{)}{0pt}{}{a+1}{b+c+1}$. 
\end{s}\begin{proof}
Denoting the left-hand side by $F(a,b,c)$ and using the
equality $\genfrac{(}{)}{0pt}{}{x}{c}=
\sum_{y=c}^x\genfrac{(}{)}{0pt}{}{y-1}{c-1}$, we obtain
\[
 F(a,b,c)=\sum_{x=c}^{a-b}\genfrac{(}{)}{0pt}{}{a-x}{b}
\sum_{y=c}^x\genfrac{(}{)}{0pt}{}{y-1}{c-1}=
\sum_{y=c}^{a-b}\genfrac{(}{)}{0pt}{}{y-1}{c-1}
\sum_{x=y}^{a-b}\genfrac{(}{)}{0pt}{}{a-x}{b}=\ldots=F(a,b+1,c-1)\ .
\]
Hence, it suffices to treat the case $c=0$, which is easy.
\end{proof}%
For $\spn$, the automorphism group of $\fH(\Ab(\spn))$ is not as rich
as for $\slno$. But a direct computation is not difficult.

\begin{s}{Theorem}  \label{edges:sp}
If $\g=\spn$, $n\ge 1$, then $\#\EAb_i=2^{n-2}$ for $i\in\{1,\ldots,n-1\}$
and $\#\EAb_0=2^{n-1}$.
\end{s}\begin{proof*}
In $\spn$, the commutative roots are 
$\gamma_{ij}:=\esi_i+\esi_j$, $1\le i\le j \le n$ and
$cl(\gamma_{ij})=\ap_{j-i}$. This already shows that there is no commutative roots
of class $\ap_{n}$, and thereby no edges of type $\ap_{n}$ in $\EAb$.

The standard matrix realisation of $\spn$ shows that the
set of commutative roots is identified with the left-justified Ferrers diagram
with row lengths $(n,n-1,\ldots,1)$, and 
the Abelian upper ideals are identified with  shifted
Ferrers subdiagrams inside, cf. also \cite[(5.1)]{duality}. 
For instance, the 7-dimensional ideal in the right diagram in Figure~\ref{ab:sp}
has the generators $\gamma_{14}$ and $\gamma_{33}$ of classes $\ap_0$ and $\ap_3$,
respectively.
For brevity, we will say that an Abelian ideal in $\spn$
is depicted by a {\it triangular\/} Ferrers diagram of size $\le n$.
\begin{figure}[htb]
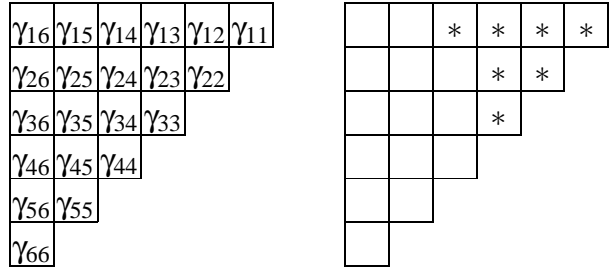

\begin{center}
\setlength{\unitlength}{0.02in}
\Yboxdim3ex  
{\young(\zas\zar\zaq\zap\zan\zam,\zbs\zbr\zbq\zbp\zbn,\zcs\zcr\zcq\zcp,\zds\zdr\zdq,\zes\zer,\zfs)} \qquad
{\young(\ \ \ast \ast \ast \ast ,\ \ \ \ast \ast ,\ \ \ \ast ,\ \ \ ,\ \ ,\ )}
\caption{The commutative roots and an Abelian ideal
for ${\frak sp}_{12}$} \label{ab:sp}
\end{center}
\end{figure}
Take the commutative root of class $\ap_i$ lying in the $k$-th row of the diagram,
i.e., $\gamma_{i,i+k}$. Then $k\le n-i$.
Using these pictures, we find that
the number of triangular Ferrers diagram of size $\le n$ having this root as
a generator (=\,south-west corner) equals
\\[.8ex]  \hspace*{2ex}
$\left\{\parbox{200pt}{the number of the (usual) Ferrers diagrams that fit in 
the rectangle of the shape $(n-i-k)\times (k-1)$}\right\}\ \times \   
\left\{\parbox{185pt}{the number of triangular Ferrers diagrams 
of size $\le i{-}1$}\right\}$. 
\\[.8ex]
Here the first number equals $\genfrac{(}{)}{0pt}{}{n-i-1}{k-1}$, and the second
one equals $2^{i-1}$, if $i\ge 1$.
Therefore the total number of 
Abelian ideals having
a generator of class $\ap_i$ is equal to
\[
  \sum_{k=1}^{n-i}\genfrac{(}{)}{0pt}{}{n-i-1}{k-1}2^{i-1}=2^{n-i-1}2^{i-1}=2^{n-2} \ .  
\]
Similarly, for $i=0$ we obtain
$\displaystyle
  \sum_{k=1}^{n}\genfrac{(}{)}{0pt}{}{n-1}{k-1}=2^{n-1} \ .  
$
\end{proof*}%
\begin{s}{Theorem}  \label{edges:soo}
If\/ $\g=\sono$ or\/ $\sone$, $n\ge 3$, then\/ $\#\EAb_i=2^{n-2}$ for each
$i\in\{0,1,\ldots,n\}$.
\end{s}\begin{proof}
1. $\g=\sono$. The classes of the commutative roots are as follows: \\
$cl(\esi_1-\esi_2)=\ap_0$,
$cl(\esi_1-\esi_j)=\ap_{j-1}$, where $3\le j\le n+1$ and $\esi_{n+1}:=0$,
and 
\begin{equation}  \label{i+j}
\mbox{$cl(\esi_i+\esi_j)=\left\{\begin{array}{cl}
\ap_{j-i}, & \text{if $j-i\ge 2$,} \\
\ap_0,     & \text{if $j-i=1$ and $j$ is even,}  \\
\ap_1,     & \text{if $j-i=1$ and $j$ is odd.}  \\
\end{array}\right.$ \ \ where $1\le i<j \le n$.}
\end{equation}
Notice that $\esi_1$ is the only short commutative roots and hence the only root
of class $\ap_n$.
The standard matrix realisation of $\sono$ shows that the
set (upper ideal) of commutative roots is identified with a skew Ferrers diagram, 
with row lengths $(2n-1,n-2,n-3,\ldots,1)$, see the sample Figure for ${\frak so}_{13}$,
where the leftmost (resp. rightmost) box corresponds to $\ap_1=\esi_1-\esi_2$
(resp. $\theta=\esi_1+\esi_2$) and the empty boxes represent the non-commutative roots.

\begin{figure}[htb]
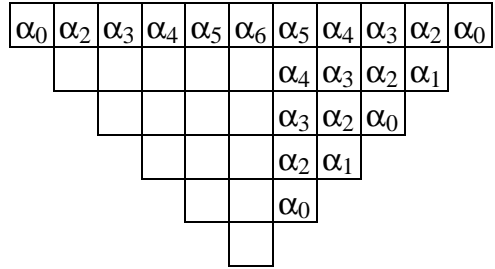

\begin{center}
\setlength{\unitlength}{0.02in}
\Yboxdim3ex  
{\young(\apo\apxx\apxxx\apxxxx\apxxxxx\apxxxxxx\apxxxxx\apxxxx\apxxx\apxx\apo ,:\ \ \ \ \ \apxxxx\apxxx\apxx\apx ,::\ \ \ \ \apxxx\apxx\apo ,:::\ \ \ \apxx\apx ,::::\ \ \apo,:::::\ )}
\caption{The classes of commutative roots 
for ${\frak so}_{13}$} \label{ab:sono}
\end{center}
\end{figure}

\noindent 
For each commutative root, we can compute the number of Abelian ideals having this 
root as a generator. (It is important that the ideal of the commutative roots is
not Abelian, so that one should count sub-ideals that are really Abelian.) 
For instance, if $i\ge 2$ and $j-i\ge 2$, then the number of Abelian
ideals having $\esi_i+\esi_j$ 
as a generator is equal to $\displaystyle 2^{j-i-1}\genfrac{(}{)}{0pt}{}{n-j+i-1}{i-1}$;
if $i\ge 2$ and $j-i=1$, then the number of ideals having $\esi_i+\esi_{i+1}$ 
as a generator is equal to $\displaystyle 2 \genfrac{(}{)}{0pt}{}{n-2}{i-1}$. 
For the roots $\esi_1\pm \esi_j$, the computations and the answer are even easier.
(We omit explicit manipulations with Ferrers diagrams.)

\noindent
Then taking the sum of the numbers obtained
corresponding to the roots of the same class yields the answer.

2. $\g=\sone$.  The classes of the commutative roots are as follows: \\
$cl(\esi_1-\esi_2)=\ap_0$,
$cl(\esi_1-\esi_j)=\ap_{j-1}$ \ ($3\le j\le n-1$), 
$cl(\esi_i-\esi_n)=\ap_{n-i}$ \ ($1\le i\le n-2$),
$cl(\esi_{n-1}-\esi_n)=cl(\esi_{n-1}+\esi_n)$, and also Eq.~\re{i+j}.

\noindent
Again, it is easier to understand this by looking at the sample picture for
${\frak so}_{12}$, which is associated with the standard matrix
realisation of $\sone$.

\begin{figure}[htb]
\begin{center}
\setlength{\unitlength}{0.02in}
\Yboxdim3ex  
{\young(\apo\apxx\apxxx\apxxxx\apxxxxxx\apxxxxx\apxxxx\apxxx\apxx\apo ,:\ \ \ \apxxxx\apxxxx\apxxx\apxx\apx ,::\ \ \apxxx\apxxx\apxx\apo ,:::\ \apxx\apxx\apx ,::::\apo\apo )}
\caption{The classes of commutative roots 
for ${\frak so}_{12}$} \label{ab:sone}
\end{center}
\end{figure}

\noindent 
We omit all other details for this case.
\end{proof}%
Theorems~\ref{edges:sl}, \ref{edges:sp}, \ref{edges:soo}, and our computations 
for the exceptional Lie algebras prove the following general result.

\begin{s}{Theorem} \label{thm:edges} 
For each simple Lie algebra $\g$, the numbers $\#\EAb_i$ depend only on the
length of $\ap_i$. Moreover,
if $\g\ne \spn$, then the number of edges of each type
in $\fH(\Ab)$ is equal to $2^{n-2}$ (\text{\normalfont Warning}: $\tri$ is regarded as ${\frak sp}_2$).
\end{s}%
Our proof for the exceptional Lie algebras
uses the same strategy together with the explicit knowledge of classes of
the commutative roots. 
The classes of all commutative roots for $\GR{E}{6},\,\GR{E}{7},\,\GR{E}{8},\,\GR{F}{4}$
are indicated in the Appendix. The case of $\GR{G}{2}$ is trivial.
In the figure for $\GR{F}{4}$, the leftmost node corresponds to the
zero-dimensional ideal and the two sinks on the right
represent two maximal Abelian ideals.  

\begin{figure}[htb]
\setlength{\unitlength}{0.03in}
\begin{center}
\begin{picture}(110,70)(1,1)

\multiput(0,40)(15,0){5}{\circle{2}}
\multiput(70,30)(10,10){4}{\circle{2}}
\multiput(80,20)(10,10){3}{\circle{2}}
\multiput(90,10)(10,10){3}{\circle{2}}
\put(70,50){\circle{2}}
\multiput(1.5,40)(15,0){4}{\vector(1,0){12}}
\multiput(61,39)(10,-10){3}{\vector(1,-1){8}}
\multiput(71,49)(10,-10){3}{\vector(1,-1){8}}
\multiput(91,49)(10,-10){2}{\vector(1,-1){8}}

\multiput(61,41)(10,-10){4}{\vector(1,1){8}}
\multiput(81,41)(10,-10){3}{\vector(1,1){8}}
\put(91,51){\vector(1,1){8}}
\multiput(60,46)(10,-10){4}{{\small $\ap_1$}}
\multiput(80,46)(10,-10){3}{{\small $\ap_2$}}
\multiput(70,23)(10,10){3}{{\small $\ap_4$}}
\multiput(80,13)(10,10){3}{{\small $\ap_0$}}
\put(90,56){{\small $\ap_3$}}
\multiput(60,33)(10,10){2}{{\small $\ap_3$}}
\put(5,42){{\small $\ap_0$}}
\put(20,42){{\small $\ap_4$}}
\put(35,42){{\small $\ap_3$}}
\put(50,42){{\small $\ap_2$}}
\put(-2,32){{\small $\varnothing$}}

\end{picture}
\end{center}
\caption{The Hasse diagram of $\Ab(\GR{F}{4})$, with types of edges}   \label{ab_F4}
\end{figure}
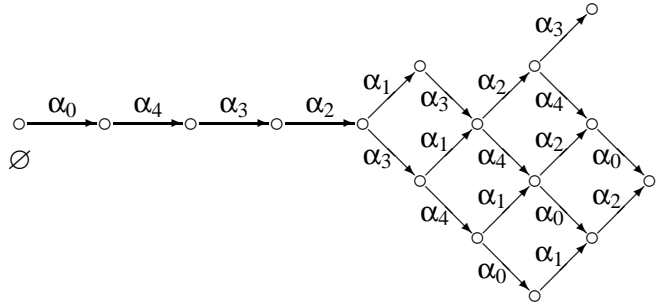
\noindent

\begin{rem}{Remark}
If $\Delta$ is not simply-laced, then one may consider 
the {\it long\/} Abelian ideals,
i.e., the Abelian ideals such that $I_\ah$ contains only long roots \cite{long}.
The corresponding poset, denoted $\Ab_l$, consists of $2,3,4$ elements for
$\GR{C}{n},\GR{G}{2},\GR{F}{4}$, respectively. The only interesting case is that 
of $\GR{B}{n}$, where $\#\Ab_l=2^{n-1}$.
It follows from \cite[Prop.\,2.1]{long} that the type of any edge of $\fH(\Ab_l)$
is a long root from $\HP$. 
One can show here that the number of edges of each `long' type is $2^{n-3}$.
In particular, $\#\mathcal E(\Ab(\GR{B}{n})_l)=n 2^{n-3}$.
\end{rem}%
\vskip-1ex


\section{The covering polynomial of a finite poset}
\label{fin_poset}
\setcounter{equation}{0}

\noindent
Results on counting edges, especially Proposition~\ref{nar+edges}, suggest a
general construction of a polynomial attached to a finite poset.

Let $(\PP, \curle)$ be a finite poset. Write $\EE(\PP)$ for the set of edges
of the Hasse diagram of $\PP$.
We define a polynomial which encodes some
properties of the covering relation in $\PP$. For any $x\in\PP$, let $\kappa(x)$
be the number of $y\in \PP$ such that $y$ is covered by $x$.

\begin{rem}{Definition}  \label{cover_pol}
The {\it covering polynomial\/} of $\PP$ is
$ \displaystyle
  \ck_{\PP}(q)=\sum_{x\in\PP} q^{\kappa(x)} \ .
$
\end{rem}%
It follows from the definition that 
\[
\ck_{\PP}(q)\vert_{q=1}=\#\PP \quad  \text{and} \quad
\frac{d}{dq}\ck_\PP(q)\vert_{q=1}=\#\EE(\PP) \ .
\]
On the other hand, $\ck_\PP(0)$ is the number of the minimal elements in $\PP$.
It might be interesting to realise whether some other values of $\ck_\PP$ have 
a combinatorial meaning.

{\bf Example.} For the Boolean lattice ${\mathcal B}_n$, we have
$\ck_{{\mathcal B}_n}(q)=(1+q)^n$.
\\[.6ex]
Let us look at the covering polynomials of posets considered above.

{\bf (A)} \ It is not hard to compute the covering polynomial for the posets $\Delta^+$.
Here is the answer:

\begin{tabular}{cl}
$\Delta$ & \phantom{quq}$\ck_{\Delta^+}(q)$ \\ \hline
$\GR{A}{n}$ & $n +\genfrac{(}{)}{0pt}{}{n}{2}q^2$ {\rule{0pt}{2.5ex}}\\
$\GR{B}{n},\GR{C}{n}$ & $n + (n-1)q+ (n-1)^2q^2$ {\rule{0pt}{2.4ex}}\\
$\GR{D}{n}$ & $n + (n-3)q+ (\genfrac{(}{)}{0pt}{}{n}{2}+\genfrac{(}{)}{0pt}{}{n-3}{2})
q^2 +(n-3)q^3$ {\rule{0pt}{2.3ex}} \\
$\GR{E}{6}$ & $6+\phantom{1} 5q+20q^2+\phantom{1} 5q^3$  {\rule{0pt}{2.3ex}}\\
$\GR{E}{7}$ & $7+10q+36q^2+10q^3$ \\
$\GR{E}{8}$ & $8+21q+70q^2+21q^3$ \\
$\GR{F}{4}$ & $4+\phantom{1} 7q+12q^2+\phantom{1} q^3$ \\
$\GR{G}{2}$ & $2+\phantom{1} 3q+\phantom{1} q^2$ \\ \hline
\end{tabular} 
\\[1ex]
Notice that {\sf (1)} $\deg \ck_{\Delta^+}\le 3$, and {\sf (2)} if 
$\Delta$ is simply-laced, then the coefficients
of $q$ and $q^3$ are equal. 
The first claim can be proved a priori; the second claim is derived from
the first one using the equalities $\ck_{\Delta^+}(1)=nh/2$ and
$\ck'_{\Delta^+}(1)=n(h-2)$, where the last equality stems from 
Theorem~\ref{simply}.

{\bf (B)} \ If $\PP=\AD(\g)$, then $\kappa(\ce)$ is the number of generators of $\ce$
for any $\ce\in\AD(\g)$. Hence the covering polynomial of 
$\AD(\g)$ is nothing but the generalised Narayana polynomial $\N_\g(q)$
from Section~\ref{pos_all}. These polynomials are explicitly written out in 
\cite[Section\,6]{duality}. See also \cite[Section\,5.2]{RW} and \cite[5.2]{bessis}, 
where $\N_\g(q)$
appears as the $h$-polynomial of the polytope $\MM(\g)$ mentioned
in Remark~\ref{politop} and a Poincar\'e polynomial for the set of simple
elements in the dual braid monoid, respectively.

{\bf (C)} \  More generally, if ${\goth X}$ is any class of \adn ideals that is closed 
under taking sub-ideals, then $\kappa(\ce)$
(for the corresponding poset) is just the number of generators of $\ce$. 
In particular, this applies to $\goth X=\Ab$.

The covering polynomial for $\PP=\Ab$ is computed for each simple Lie algebra
separately. We first note that the coefficient of $q$ is the number of
commutative roots.
It follows from a result of Sommers \cite[Theorem\,6.4]{eric} that 
the number of generators of an Abelian ideal
is at most the maximal number of orthogonal roots in $\Pi$.
This provides an upper bound on $\deg\ck_\Ab$. Actually, it is easy to
find out an Abelian ideal with such a number of generators, so that one obtains 
the precise value of the degree. For instance, $\deg\ck_{\Ab(\GR{D}{n})}=
\left[\frac{n}{2}\right]+1$,  $\deg\ck_{\Ab(\GR{E}{7})}= \deg\ck_{\Ab(\GR{E}{8})}=4$.
 
If $\deg\ck_\Ab \le 3$, then this polynomial can be computed using
the values $\ck_\Ab(1)=2^n$, $\ck'_\Ab(1)=(n+1)2^{n-2}$, and the known number of
commutative roots. 
The result for $\GR{A}{n}$ and $\GR{C}{n}$ can be deduced from matrix descriptions
of generators given in \cite[Section\,3]{pr}.
%
\begin{figure}[htb]
\begin{tabular}{cl}
$\g$ & \phantom{quq}$\ck_{\Ab}(q)$ \\ \hline
$\GR{A}{n},\GR{B}{n},\GR{C}{n}$ & $\sum_{k\ge 0}\genfrac{(}{)}{0pt}{}{n+1}{2k}q^k$ 
{\rule{0pt}{2.7ex}}\\
$\GR{D}{n}$ & 
$\sum_{k\ge 0}\bigl(\genfrac{(}{)}{0pt}{}{n+2}{2k}-
4\genfrac{(}{)}{0pt}{}{n-1}{2k-2}\bigr )q^k$ 
{\rule{0pt}{3ex}} \\
$\GR{E}{6}$ & $1+25q+\phantom{1} 27q^2+11q^3$  {\rule{0pt}{2.3ex}} \\
$\GR{E}{7}$ & $1+34q+\phantom{1} 60q^2+30q^3+\phantom{1} 3q^4$ {\rule{0pt}{2.3ex}}  \\
$\GR{E}{8}$ & $1+44q+118q^2+76q^3+17q^4$ {\rule{0pt}{2.3ex}} \\
$\GR{F}{4}$ & $1+10q+\phantom{11} 5q^2$ \\
$\GR{G}{2}$ & $1+\phantom{1} 3q$ \\
\hline
\end{tabular} 
\end{figure}
It is a mere coincidence
that the covering polynomial for $\Ab(\GR{B}{n})$ is the same
as for $\GR{A}{n}$ and $\GR{C}{n}$. The formulae for both orthogonal series 
result from direct matrix considerations, 
which are similar 
to that in Theorems~\ref{edges:sl},\,\ref{edges:soo} and \cite[Section\,3]{pr}.
The proof for $\GR{D}{n}$ is more involved, but still it can be derived
from the consideration of the Ferrers diagram for the \adn ideal of commutative roots.
(We came up with the explicit formula for $\GR{D}{n}$ using computer 
calculations made by G.\,R\"ohrle for $\GR{D}{n}$, $n=7,8$.)
The details of these computations and some related results on covering polynomials
will appear elsewhere.


\appendix
\section{Classes of the commutative roots for $\GR{E}{6}$,\,
$\GR{E}{7}$,\,$\GR{E}{8}$,\,$\GR{F}{4}$}

\noindent
We follow the natural numbering of simple roots, i.e., that adopted in
\cite{VO}. To convey the idea of "naturality" to the interested
reader, we present this graphically:

$\GR{E}{6}$:\quad 
\begin{tabular}{@{}c@{}}
1-2-\lower3.1ex\vbox{\hbox{3\rule{0ex}{2.4ex}}
\hbox{\hspace{0.5ex}\rule{.1ex}{1ex}\rule{0ex}{1.2ex}}\hbox{6\strut}}-4-5
\end{tabular} 
\qquad
$\GR{E}{7}$:\quad 
\begin{tabular}{@{}c@{}}
1-2-3-\lower3.1ex\vbox{\hbox{4\rule{0ex}{2.4ex}}
\hbox{\hspace{0.5ex}\rule{.1ex}{1ex}\rule{0ex}{1.2ex}}\hbox{7\strut}}-5-6
\end{tabular} 
\qquad
$\GR{E}{8}$:\quad 
\begin{tabular}{@{}c@{}}
1-2-3-4-\lower3.1ex\vbox{\hbox{5\rule{0ex}{2.4ex}}
\hbox{\hspace{0.5ex}\rule{.1ex}{1ex}\rule{0ex}{1.2ex}}\hbox{8\strut}}-6-7
\end{tabular}
\qquad
$\GR{F}{4}$:\quad  1-2$\Leftarrow$3-4 .

\noindent
For brevity, we write $[m_1\,m_2\ldots m_n]$ in place of
$\sum_{i=1}^n m_i\ap_i$.

\begin{itemize}
\item[$\GR{E}{6}$:]
\begin{itemize}
\item[$\ap_0$] - [111110], [123212];
\item[$\ap_1$] - [000010], [012211]; 
\item[$\ap_2$] - [000110], [111001], [012111], [112211];
\item[$\ap_3$] - [111000], [001110], [111101], [011111], [012101], [112111], [122211];
\item[$\ap_4$] - [110000], [001111], [112101], [122111];
\item[$\ap_5$] - [100000], [122101];
\item[$\ap_6$] - [111100], [011110], [111111], [123211].  
\end{itemize}
\item[$\GR{E}{7}$:]
\begin{itemize}
\item[$\ap_0$] - [1000000], [1222101], [1234322]; \rule{0pt}{3ex}
\item[$\ap_1$] - [0123212];
\item[$\ap_2$] - [1111110], [0123211], [1123212];
\item[$\ap_3$] - [1111100], [1111111], [0122211], [1123211]. [1223212];
\item[$\ap_4$] - \parbox[t]{390pt}{[1111000], [1111101], [1112111], [0122111], [0112211],
                \rule{0pt}{2.3ex}[1122211],[1223211], [1233212];}
\item[$\ap_5$] - [1110000], [1112101], [0122101], [0012211], [1122111], [1222211],
                 [1234212]\rule{0pt}{2.3ex}; 
\item[$\ap_6$] - [1100000], [1122101], [1222111], [1234312];
\item[$\ap_7$] - [1111001], [1112211], [1233211].
\end{itemize}
\item[$\GR{E}{8}$:]
\begin{itemize}
\item[$\ap_0$] - [01234322], [23456423]; \rule{0pt}{3ex}
\item[$\ap_1$] - [01234312], [11234322], [13456423];
\item[$\ap_2$] - [12222101], [01234212], [11234312], [12234322], [12456423];
\item[$\ap_3$] - [12222111], [01233212], [11234212], [12234312], [12334322], [12356423];
\item[$\ap_4$] - \parbox[t]{370pt}{[01233211], [01223212], [12222211], [11233212], [12234212], [12334312],\rule{0pt}{2.3ex}
                 [12344322], [12346423];}
\item[$\ap_5$] - \parbox[t]{370pt}{[01123212], [11233211], [11223212], [12223211], [12233212], [12334212],\rule{0pt}{2.3ex}
                 [12344312], [12345322], [12345423];} 
\item[$\ap_6$] - [11123212], [12233211], [12333212], [12345312], [12345323]
\rule{0pt}{2.3ex};
\item[$\ap_7$] - [12333211], [12345313];
\item[$\ap_8$] - [00123212], [12223212], [12344212], [12345422].
\end{itemize}
\item[$\GR{F}{4}$:]
\begin{itemize}
\item[$\ap_0$] - [2210], [2432]; \rule{0pt}{3ex}
\item[$\ap_1$] - [1321];
\item[$\ap_2$] - [1221], [2321]; 
\item[$\ap_3$] - [0221], [2221], [2421];
\item[$\ap_4$] - [2211], [2431].
\end{itemize}
\end{itemize}
To obtain these lists, we combine explicit computations
of certain minimal elements with using  Proposition~\ref{reduced}
and Lemma~\ref{sosedi}.
The information for $\GR{F}{4}$ is also readily being extracted from the data 
of Table~6.1 in \cite{lp}.

\end{document}